\documentclass[a4paper]{article}

\usepackage[utf8]{inputenc} \usepackage{lmodern} \usepackage[ngerman, english]{babel} \usepackage{csquotes} \usepackage{mathrsfs}

\usepackage[unicode, colorlinks]{hyperref} \usepackage{enumitem}

\usepackage{xcolor}

\usepackage{booktabs}

\usepackage{mathtools} \usepackage{amsthm, amssymb} \usepackage{aligned-overset}

\usepackage{url}
\usepackage[numbers]{natbib}

\usepackage{comment}
 \usepackage[left=3.6cm,right=3.6cm,bottom=3cm,top=3cm]{geometry}
\usepackage{orcidlink}

\theoremstyle{plain}\newtheorem{prop}{Proposition}[section]
\newtheorem{lemma}[prop]{Lemma}
\newtheorem{corollary}[prop]{Corollary}
\newtheorem{theorem}[prop]{Theorem}

\theoremstyle{definition}
\newtheorem{defenv}[prop]{Definition}

\newtheorem{remenv}[prop]{Remark}
\newtheorem{exenv}[prop]{Example}

\newenvironment{remark}
{\pushQED{\qed}\remenv}
{\popQED\endremenv}
\newenvironment{definition}
{\pushQED{\qed}\defenv}
{\popQED\endremenv}
\newenvironment{example}
{\pushQED{\qed}\exenv}
{\popQED\endremenv}

\newcommand*{\dims}{d}

\newcommand*{\real}{\mathbb{R}}
\newcommand*{\nat}{\mathbb{N}}
\newcommand*{\integer}{\mathbb{Z}}
\newcommand*{\complex}{\mathbb{C}}

\newcommand*{\sphere}[1][\dims-1]{\mathbb{S}^{#1}}
\newcommand*{\domain}{\mathbb{X}}

\newcommand{\E}{\mathbb{E}}
\renewcommand{\Pr}{\mathbb{P}}

\DeclareMathOperator{\Cov}{Cov}
\newcommand*{\normal}{\mathcal{N}}

\DeclareMathOperator{\iid}{iid}

\newcommand*{\kernel}{K}
\newcommand*{\PDset}{\mathcal{P}}

\newcommand*{\rf}{\mathbf{f}}
\newcommand*{\rg}{\mathbf{g}}

\newcommand{\identity}{\mathrm{Id}}
\newcommand*{\transpose}{T}

\newcommand*{\conj}[1]{\overline{#1}}

\DeclarePairedDelimiterXPP{\scp}[2]{}{\langle}{\rangle}{}{#1, #2}
\DeclarePairedDelimiter{\norm}{\|}{\|}
\DeclarePairedDelimiter{\abs}{|}{|}
\DeclarePairedDelimiter{\set}{\{}{\}}
\DeclarePairedDelimiter{\floor}{\lfloor}{\rfloor}

\DeclareMathOperator{\arcosh}{arcosh}
\newcommand*{\metric}{\rho}

\newcommand*{\uP}[1][\lambda]{C^{#1}} \newcommand*{\uPnorm}[1][\lambda]{\hat{C}^{#1}} \newcommand*{\weight}{w}

\definecolor{primary}{RGB}{0,48,86}

\colorlet{primary10}{primary!10}
\colorlet{primary25}{primary!25}
\colorlet{primary40}{primary!40}
\colorlet{primary55}{primary!55}
\colorlet{primary70}{primary!70}
\colorlet{primary85}{primary!85}

\definecolor{accent}{RGB}{179,182,185}

\newcommand*{\red}[1]{{\color{red} #1}}

\newcommand{\xx}{{x}}
\newcommand{\xy}{{y}}

\newcommand{\xxs}{\xx}\newcommand{\xys}{y}
\newcommand{\nngp}[1]{\kernel^{\mathrm{NNGP}}_{#1}}
\newcommand{\ntk}[1]{\kernel^{\mathrm{NTK}}_{#1}}
\newcommand{\pdkernel}{positive definite kernel}
\newcommand{\pdkernels}{positive definite kernels}
\newcommand{\shortkernel}{kernel}
\newcommand{\shortkernels}{kernels}

\newcommand{\spdkernels}{strictly positive definite kernels}

\newcommand{\equalDef}{\coloneqq}
\newcommand{\defEqual}{\eqqcolon}

\newcommand*{\homP}[1][\dims]{\mathcal{P}^{#1}}
\newcommand*{\homPSphere}[1][\dims]{\mathscr{P}^{#1}}
\newcommand*{\homHarmonic}[1][\dims]{\mathcal{H}^{#1}}
\newcommand*{\sphHarmonic}[1][\dims]{\mathscr{H}^{#1}}
\newcommand*{\zonalK}[1][\dims]{Z^\dims}
\DeclarePairedDelimiterXPP{\vol}[1]{}{|}{|}{}{#1}

\renewcommand*{\restriction}{\mathord{\upharpoonright}}
\newcommand*{\besselProj}{\mathcal{J}}
\newcommand*{\bessel}{J}

\title{Schoenberg characterization of continuous non-stationary isotropic positive definite kernels}
\author{
	Felix Benning\thanks{Corresponding author}
	~\orcidlink{0009-0001-2354-7696}\\
	\texttt{felix.benning@gmail.com}\\
	University of Mannheim
	\and
	Max David Schölpple
	\orcidlink{0009-0000-9078-4205}
	\\
	\texttt{max.schoelpple@mathematik.uni-stuttgart.de}\\
	University of Stuttgart
}

\begin{document}
	\maketitle 

	\begin{abstract}
	We characterize the continuous isotropic positive definite kernels on \(\real^\dims\),
	where isotropy refers to invariance under the orthogonal
	group \(O(\dims)\) but not necessarily stationarity. Furthermore, we
	characterize strict positive definiteness for such kernels.

	\noindent The class of isotropic kernels is fairly general as it unifies
	stationary isotropic and dot product kernels, and includes
	neural network kernels that arise from infinite-width limits of neural
	networks.
	As an application, we further characterize the continuous isotropic Gaussian
	random functions in terms of a series representation.
	\smallskip

	\noindent \textbf{Keywords}: positive definite, kernels, isotropy, characterization,
	non-stationary, neural network kernels

	\noindent \textbf{2020 MSC}: 33C50, 33C55, 42A82, 42C10, 43A35, 60G15, 68T07
\end{abstract} 	\section{Introduction}
\label{sec: introduction}

We call a function \(\kernel\colon \domain \times \domain \to \complex\) a \emph{\pdkernel}, or short \emph{\shortkernel}, if for all \(m\in \nat\),  \(x_1,\dots, x_m\in \domain\) and \(c=(c_1,\dots,c_m) \in \complex^m\)
\begin{align}
	    \sum_{i,j=1}^m c_i \conj{c_j}\kernel(x_i, x_j) \ge 0.
	\label{eq:positive_definite_definition}
\end{align}
If the inequality in \eqref{eq:positive_definite_definition} is strict whenever $x_1,\dots,x_m$ are distinct vectors and \(c \neq 0\), we call the \shortkernel{}
\emph{strictly positive definite}. Recall that any
\shortkernel{} that satisfies \eqref{eq:positive_definite_definition} is
Hermitian, that is $\kernel(x,y) = \conj{\kernel(y,x)}$ \citep[Lem.~3]{berlinetReproducingKernelHilbert2004}.  Numerous previous works
established characterizations for classes of \shortkernels{} invariant to
certain input transformations \citep[e.g.][]{bochnerMonotoneFunktionenStieltjessche1933,
	schoenbergMetricSpacesPositive1938,schoenbergPositiveDefiniteFunctions1942,guellaSchoenbergsTheoremPositive2019,
	estradeCovarianceFunctionsSpheres2019,bergSchoenbergCoefficientsSchoenberg2017}. 
These transformations are typically groups
acting on the space \(\domain\). In the case \(\domain =\real^\dims\)
the most natural groups to consider are perhaps
\begin{itemize}[noitemsep]
    \item the group of translations \(T(\dims) := \{\phi\colon x\mapsto x+ v \mid v \in \real^\dims\}\),
    \item the group of Euclidean  isometries
    \[
        E(\dims)
        := \bigl\{\phi\colon \real^\dims \to \real^\dims \mid \forall x,y\in \real^\dims, \; \|\phi(x)-\phi(y)\| = \|x-y\|\bigr\},
    \]
    \item and the group of linear isometries \(O(\dims):= \{\phi \in E(\dims) : 
    \phi(0)=0\}\), which is canonically represented by the group of orthogonal
    matrices \(U\) that satisfy $UU^\transpose=\identity$.
\end{itemize}
The group of linear isometries $O(d)$ readily generalizes to the infinite
dimensional space of square summable sequences
\begin{equation}
    \label{eq: definition of l2}    
    \real^\infty \equalDef \ell^2
    \equalDef \Bigl\{
        (x_k)_{k\in \nat}\in \real^\nat \Bigm| \|x\|^2 := \sum_{k\in \nat} x_k^2< \infty\Bigr\}
\end{equation}
and we denote this group by \(O(\infty)\). 
Similarly, we denote the group of Euclidean isometries on \(\real^\infty\) by
\(E(\infty)\).
Using such groups we may define a general notion of invariant kernels.

\begin{definition}[\(\Phi\)-invariance]
	Let $\Phi$ be a set of transformations $\phi:\domain\to\domain$, typically
    a group. A \pdkernel{} \(\kernel\) is called \(\Phi\)-invariant if
    \[
        \kernel(\phi(x), \phi(y)) = \kernel(x,y) \qquad \forall \phi \in \Phi,\; \forall x,y\in \domain,
    \]
    and we denote the set of \(\Phi\)-invariant \shortkernels{} by \(\PDset_\Phi:=\PDset_\Phi(\domain)\).
    In the case \(\domain = \real^\dims\) we say the \shortkernel{} \(\kernel\) is
    \begin{itemize}[noitemsep]
        \item \emph{stationary}, if \(\Phi\) is the group of translations
        \(T(\dims)\),
\item \emph{stationary isotropic}, if \(\Phi\) is the group of (Euclidean) isometries \(E(\dims)\).
        \item \emph{isotropic}, if \(\Phi\) is the orthogonal group \(O(\dims)\) of linear isometries. \qedhere
    \end{itemize}
\end{definition}

A \shortkernel{} which is invariant to a set of transformations \(\Phi\)
is also invariant to any subset of \(\Phi\), that is \(\Psi
\subseteq \Phi\) implies \(\PDset_\Psi
\supseteq \PDset_\Phi\). 
The set of stationary isotropic \shortkernels{}
\(\PDset_{E(\dims)}\) is consequently a subset of \(\PDset_{T(\dims)}\) and \(\PDset_{O(\dims)}\).
While there have long been characterizations for the stationary \shortkernel{} \(\PDset_{T(\dims)}\)
\citep{bochnerMonotoneFunktionenStieltjessche1933} and the stationary isotropic
\shortkernels{} \(\PDset_{E(\dims)}\) \citep{schoenbergMetricSpacesPositive1938}, the
non-stationary isotropic \shortkernels{} \(\PDset_{O(\dims)}\) are, to the best of our knowledge, 
left without characterization. In this work we fill this gap. 

\begin{table*}[h]
	\caption[]{
		Characterizations of continuous \pdkernels{} \(\kernel\colon \domain\times \domain \to \complex\).
		We use \(\Omega_\dims(t):=\E[e^{itU_1}], U\sim \mathrm{Unif}(\sphere)\) 
and \(\uPnorm_n\) denotes the normalized
		Gegenbauer polynomial of degree \(n\) with parameter \(\lambda =
		\frac{\dims-2}{2}\), cf.~\eqref{eq: Gegenbauer polynomials}.
	}
	\label{table: characterizations}
	\centering
	\makebox[\textwidth][c]{
		\(\def\arraystretch{1.2}
		\begin{array}{l l l l l}
			\domain & \Phi & \kernel(x,y) & \text{Characterization}
			\\
			\toprule
			\real^\dims
			& T(\dims)
			& \kappa(x-y)
			& \kappa(u) = \int e^{i \langle v, u\rangle} d\mu(v)
			& \mu \text{ finite measure}
			\\
			\real^\dims
			& E(\dims)
			& \kappa(\|x-y\|)
			& \kappa(r) = \int_0^\infty \Omega_\dims(rs) d\mu(s)
			& \mu \text{ finite measure}
			\\
			\ell^2
			& E(\infty)
			& \kappa(\|x-y\|)
			& \kappa(r) = \int_0^\infty \exp(-s^2r^2/2) d\mu(s)
			& \mu \text{ finite measure}
			\\
			\sphere
			& O(\dims)
			& \kappa(x\cdot y)
			& \kappa(\rho) = \sum_{n=0}^\infty \gamma_n \uPnorm_n(\rho)
			& \gamma_n\ge 0
			\\
			\sphere[\infty]
			& O(\infty)
			& \kappa(x\cdot y)
			& \kappa(\rho) = \sum_{n=0}^\infty \gamma_n \rho^n
			& \gamma_n\ge 0
			\\
			\mathcal X \times \sphere
			& \identity_{\mathcal X} \otimes O(\dims)
			& \kappa(x_{\mathcal X}, y_{\mathcal X}, x_{\sphere[]} \cdot y_{\sphere[]})
			& \kappa(r,s,\rho) = \sum_{n=0}^\infty \alpha_n(r,s) \uPnorm_n(\rho)
			& \alpha_n \text{ pos. definite}
			\\
			\mathcal X \times \sphere[\infty]
			& \identity_{\mathcal X} \otimes O(\infty)
			& \kappa(x_{\mathcal X}, y_{\mathcal X}, x_{\sphere[]} \cdot y_{\sphere[]})
			& \kappa(r,s,\rho) = \sum_{n=0}^\infty \alpha_n(r,s) \rho^n
			& \alpha_n \text{ pos. definite}
			\\
			\real^\dims \times \sphere
			& T(\dims) \otimes O(\dims)
			& \kappa(x_\real - y_\real, x_{\sphere[]} \cdot y_{\sphere[]})
			& \kappa(r,s,\rho) = \sum_{n=0}^\infty \alpha_n(r-s) \uPnorm_n(\rho)
			& \alpha_n \text{ pos. definite}
			\\
			\real^\dims \times \sphere[\infty]
			& T(\dims) \otimes O(\infty)
			& \kappa(x_\real - y_\real, x_{\sphere[]} \cdot y_{\sphere[]})
			& \kappa(r,s,\rho) = \sum_{n=0}^\infty \alpha_n(r-s) \rho^n 
			& \alpha_n \text{ pos. definite}
			\\
			\midrule
			\real^\dims
			& O(\dims)
			& \kappa(\|x\|, \|y\|, x\cdot y)
			& \kappa(r,s, \rho) = \sum_{n=0}^\infty \alpha_n(r,s) \uPnorm_n(\tfrac{\rho}{rs})
			& \substack{
				\text{\(\alpha_n\) pos. definite \&}\\
				\text{zero at origin}
			}
			\\
			\real^\infty
			& O(\infty)
			& \kappa(\|x\|, \|y\|, x\cdot y)
			& \kappa(r,s, \rho) = \sum_{n=0}^\infty \alpha_n(r,s)\bigl(\tfrac{\rho}{rs}\bigr)^n 
			& \substack{
				\text{\(\alpha_n\) pos. definite \&}\\
				\text{zero at origin}
			}
		\end{array}
		\)
	}
\end{table*}

To contextualize our results, Table \ref{table: characterizations} summarizes
closely related characterizations, including our own characterization of
isotropic \shortkernels{} in Theorem \ref{thm: isotropy characterization}.
The first three results can be found in the book by
\citet{sasvariMultivariateCharacteristicCorrelation2013}. The very first result
goes back to
\citet{bochnerMonotoneFunktionenStieltjessche1933} and can be generalized to
groups. Results 2-5 go back to
\citet{schoenbergMetricSpacesPositive1938,schoenbergPositiveDefiniteFunctions1942}.
In his honor results 6-9 have been called Schoenberg characterizations \citep{guellaSchoenbergsTheoremPositive2019, estradeCovarianceFunctionsSpheres2019,bergSchoenbergCoefficientsSchoenberg2017}.
The last two characterizations are our own contribution.
Many of the characterizations in Table \ref{table: characterizations}, in
particular our own, are simpler in the case \(\dims=\infty\), which describes kernels  on the Hilbert space of square summable sequences $\real^\infty := \ell^2$. 
\citet{schoenbergPositiveDefiniteFunctions1942} already observed that the kernels on \(\ell^2\) may
be viewed as a well-behaved subset of isotropic kernels in finite dimension. 
The following remark recapitulates this notion in our setting for \(\PDset_{O(\dims)}\). These ideas, formalized in Remark \ref{rem:dimensional_hierarchy} and
Lemma \ref{lem: valid in all dimensions}, readily generalize to the other natural groups
of transformations discussed in Table \ref{table: characterizations}.

\begin{remark}[Embedding]\label{rem:dimensional_hierarchy}
    Embedding the finite dimensional spaces \(\real^\dims\) into the space of
    sequences \(\real^\infty:=\ell^2\) via \(\real^\dims \simeq \set{(x_n)_{n\in \nat} \in \ell^2 : x_n = 0 \text{ for } n > \dims}\)
	naturally results in
	\begin{equation}
		\label{eq: dimensional hierarchy}	
		\real^1 \subseteq \real^2 \subseteq \dots \subseteq \real^\infty.
	\end{equation}
	An isotropic kernel $\kernel\in \PDset_{O(\dims)}$ may be identified with a
	function \(\kappa\colon \real^3 \to \complex\) such that \(\kernel(x,y) = \kappa(\|x\|, \|y\|, \langle
	x,y \rangle)\). This \(\kappa\) may be used to lift the kernel
	into any $\real^m$ for
	$m\in\nat\cup\{\infty\}$.
	The resulting bivariate function is however not necessarily positive definite
	for $m > \dims$. 
	However, \eqref{eq: dimensional hierarchy} yields $\kernel\in \PDset_{O(n)}$ for any $n\le \dims$ resulting in
    \begin{equation}
        \label{eq: pd subset relation}    
        \PDset_{O(1)}(\real) \supseteq \PDset_{O(2)}(\real^2) \supseteq \dots \supseteq \PDset_{O(\infty)}(\real^\infty).
		\qedhere
    \end{equation}
\end{remark}
The set of functions $\kappa$ that result in \emph{positive definite} kernels $\kernel(x,y)$ in all
finite dimensions is the set of isotropic kernels on the space of
eventually zero sequences
\(\bigcup_{\dims\in \nat} \real^\dims\). The following lemma shows that it coincides with the set of isotropic kernels on $\ell_2$.
\begin{lemma}
	\label{lem: valid in all dimensions}
	If $\kernel(x,y) \coloneq \kappa(\norm x,\norm y,\scp xy)$ defines a \pdkernel{} on $\real^\dims$ for all $\dims\in\nat$, then it also defines a \pdkernel{} on $\real^\infty = \ell^2$, that is
	\[
	\textstyle
	\bigcap_{\dims \in \nat} \PDset_{O(d)}(\real^\dims) = \PDset_{O(\infty)} (\ell^2).
	\]
\end{lemma}
A proof can be found in Section \ref{sec: proofs}.

\paragraph*{Outline}
In Section \ref{sec:positive_definite}, we provide the formal statement of our
characterization of continuous isotropic \shortkernels\ \(\PDset_{O(\dims)}\)
as listed in Table \ref{table: characterizations}. This characterization
essentially extends a characterization of kernels on the space \((0, \infty) \times \sphere
\simeq \real^\dims \setminus \set{0}\) by
\citet{guellaSchoenbergsTheoremPositive2019} to the origin -- as
does Section~\ref{sec: strict positive definite}, where we provide
necessary and sufficient conditions
for isotropic \shortkernels{} to be strictly positive definite.
In Section~\ref{sec: applications} we connect our results
to well known examples of kernels.
In particular we
demonstrate how the stationary isotropic kernels and dot product kernels are
embedded into our characterization and examine `neural network kernels'
arising from the theoretical analysis of neural networks.
These kernels are isotropic yet non-stationary and do not fit into any of the
previous kernel classifications. They however naturally fall into the class of
isotropic kernels considered in this work. We discuss a few more examples of
non-stationary isotropic kernels, including a
hyperbolic distance kernel, before we proceed to a first application
of our characterization in Section~\ref{sec: GRF representation}.
There we prove a series representation of isotropic Gaussian random functions
using our characterization of isotropic kernels.
We prove our main results in Section~\ref{sec: proofs}. 
In Section~\ref{sec: spherical harmonics} we recapitulate
foundations of spherical harmonics and Gegenbauer polynomials and
translate between different conventions used in the literature.

 	\section{Characterization of isotropic kernels on \texorpdfstring{\(\real^\dims\)}{Rd} }
\label{sec:positive_definite}

The defining property of isotropic \shortkernels{} is the identity
\(\kernel(Ux,Uy) = \kernel(x,y)\) for all \(x,y\in \real^\dims\) and \(U\in
O(\dims)\). The tuple \((x,y)\) may thus be characterized by the invariants
\[
	(\norm x,\norm y,\scp xy)
	\quad \text{or alternatively} \quad
	\bigl(\norm x,\norm y,\scp{\tfrac{x}{\|x\|}}{\tfrac{y}{\|y\|}}\bigr),
\]
leading to natural equivalence classes of tuples for isotropy.
In our main characterization result (Theorem \ref{thm: isotropy
characterization}) the third, spherical component
translates to Gegenbauer polynomials.
Gegenbauer polynomials \citep[e.g.][]{reimerGegenbauerPolynomials2003}
were already used by
\citet{schoenbergPositiveDefiniteFunctions1942} to characterize the
\(O(\dims)\)-invariant \shortkernels{} on \(\sphere\) and
we recapitulate their relation to spherical harmonics in Section \ref{sec: spherical harmonics}.
We denote the Gegenbauer polynomial of index \(\lambda = (\dims-2)/2\) and
degree \(n\) by \(\uP_n\). In our characterizations we use the
\emph{normalized Gegenbauer polynomials} defined by\footnote{
    Note that \(|\uP_n(t)|\le \uP_n(1)\) and \(0<\uP_n(1)=\binom{2\lambda+n-1}{n} < \infty\) for \(\lambda \neq 0\) \citep[][Eq.
    (2.16), (2.13)]{reimerGegenbauerPolynomials2003}.
}
\begin{equation}
    \label{eq: Gegenbauer polynomials}
    \uPnorm_n(t) := \begin{cases}
        \frac{1}{\uP_n(1)}\uP_n(t) & \lambda < \infty
        \\
        t^n & \lambda = \infty,
    \end{cases}
\end{equation}
where we follow the convention, that for \(\dims=2\) (\(\lambda=0\)), the Gegenbauer
polynomials are the Chebyshev polynomials
\citep[e.g.][(2.8)]{reimerGegenbauerPolynomials2003} up to normalization.

\begin{remark}[Interpretation of \(\lambda=\infty\)]
Note that the case \(\lambda=\infty\) can also be viewed as a limiting case.
Specifically, \citet[Lemma 1]{schoenbergPositiveDefiniteFunctions1942}
shows for any $t\in [-1,1]$ that \(\uPnorm_n(t)\) converges to \(t^n\) uniformly in \(n\) as
\(\lambda\to\infty\). In our setting, the index \(\lambda=(d-2)/2\)
is determined by the dimension \(\dims\). In high dimension \(\dims \gg 1\),
the set of isotropic \shortkernels{} \(\PDset_{O(\dims)}\) is therefore well approximated by
the isotropic \shortkernels{} in the space of sequences \(\real^\infty \equiv
\ell^2\), which have a more intuitive characterization in terms of monomials in
place of Gegenbauer polynomials.
\end{remark}

Recall that the isotropic \shortkernels{}
``valid in all finite dimensions'' coincide with the isotropic \shortkernels{}
on \(\ell^2\) (Lemma~\ref{lem: valid in all dimensions}).
For the interpretation of the following characterization it is therefore
often reasonable to assume \(\dims=\infty\).

\begin{theorem}[Characterization of continuous isotropic positive definite \shortkernels]
    \label{thm: isotropy characterization}
    Let \(d \in \{2,3,\dots, \infty\}\) and \(\lambda = \tfrac{\dims-2}2\).
    For a function \(\kernel\colon \real^\dims \times \real^\dims \to
    \complex\), the following assertions are
    equivalent:
    \begin{enumerate}[label={(\roman*)}]
        \item\label{it: positive, definite isotropic}
        The function $\kernel$ is a continuous, isotropic \pdkernel{}, where isotropy means
        \[
            \kernel(x,y) = \kernel(Ux, Uy)
        \]
        for every linear isometry \(U\in O(\dims)\).
        \item\label{it: pos definite, weak representation}
        The function $\kernel$ is a \pdkernel{} and has the representation
        \[
            \kernel(x,y) = \kappa\bigl(\|x\|, \|y\|, \langle x, y\rangle\bigr)
        \]
        for a continuous \(\kappa\colon M \to \complex\) with \(M=\{ (r,s, \gamma) \in [0,\infty)^2 \times \real : |\gamma| \le rs\}\).
        \item\label{it: representation}
        The function $\kernel$ has the representation
        \begin{equation}
            \label{eq: kernel representation}
            \kernel(x,y)
            = \begin{cases}
                \sum_{n=0}^\infty \alpha^{\dims}_n \bigl(\|x\|, \|y\|\bigr)
                \uPnorm_n\bigl(\bigl\langle \tfrac{x}{\|x\|},
                \tfrac{y}{\|y\|}\bigr\rangle\bigr),
                & x\neq 0 \text{ and } y\neq 0,
                \\
                \alpha_0^\dims(\norm{x}, \norm{y})
                & x= 0 \text{ or } y= 0,
            \end{cases}
        \end{equation}
        where \(\alpha^{\dims}_n\colon [0,\infty)^2\to\complex\) are unique, continuous,
        \pdkernels, with the properties: 
        \begin{enumerate}[label=(\alph*), ref=(\roman{enumi})(\alph*)]
            \item\label{it: locally uniformly convergent}
            Their \textbf{partial sums are locally uniformly convergent}, meaning that for all \(r\in [0, \infty)\) there exists a neighborhood \(U_r\subset [0,\infty) \) of \(r\) such that
            \[
                \sup_{s\in U_r}\sum_{n=N}^\infty \alpha^{\dims}_n(s,s) \to 0 \qquad \text{for}\quad N\to\infty.
            \]
            \item\label{it: zero at the origin}
            The kernels are \textbf{zero at the origin for all \(n\ge 1\)}, 
            that is
            \begin{align}
                \label{eq: zero in origin}
                \alpha^\dims_n(r,0) = \alpha^\dims_n(0,r) = 0 \qquad \text{for all }r\in[0,\infty).
            \end{align}
            
        \end{enumerate}
\end{enumerate}
    The unique kernels $\alpha_n^d$ are real-valued if and only if $\kernel$ is real-valued.
    For \(\dims <\infty\), the \shortkernels{} \(\alpha^{\dims}_n\) can be
    represented as
    \begin{align}
    	\label{eq:alpha_formula}
        \alpha^{\dims}_n(r,s)
        = \frac{\bigl\langle \kappa(r,s, rs\,\cdot), \uPnorm_n\bigr\rangle_{\weight}}{\|\uPnorm_n(\cdot)\|_{\weight}^2}
        = \frac{\int_{-1}^1 \kappa(r,s, rs\rho) \uPnorm_n(\rho)\weight(\rho) d\rho}{\|\uPnorm_n\|_{\weight}^2} ,
    \end{align}
    with weight function \(\weight(\rho) =(1-\rho^2)^{\lambda-\frac12}\)
    and inner product defined by 
    \[
        \langle f,g\rangle_{\weight} :=  \int_{-1}^1 f(\rho) \overline{g(\rho)} \weight(\rho) d\rho .
    \]
\end{theorem}

The core idea of the proof of Theorem~\ref{thm: isotropy characterization}
(Section \ref{sec:all-proofs}) is to translate
results about \pdkernels{} on \((0,\infty) \times \sphere\), independently
discovered by \citet{guellaSchoenbergsTheoremPositive2019} and
\citet{estradeCovarianceFunctionsSpheres2019}, to a characterization of
isotropic \shortkernels{} on \(\real^\dims \setminus \{0\}\). Then we
carefully extend this characterization to \(\real^\dims\) using continuity
arguments.

\begin{remark}[Characterization of isotropic kernels on balls]
    \label{rem: isotropic kernels on balls}
    The same characterization as in Theorem \ref{thm: isotropy characterization}
    holds for isotropic kernels defined on open or closed balls \(B_r\) of
    radius \(r>0\). That is, for a continuous function \(\kernel\colon B_r
    \times B_r \to \complex\), assertions \ref{it: positive, definite isotropic}, \ref{it: pos definite, weak representation} and
    \ref{it: representation} of Theorem \ref{thm: isotropy characterization} are still equivalent, after appropriately changing domains.
For this, simply replace the map to \((0, \infty) \times \sphere\) used in the proof of
    Theorem~\ref{thm: isotropy characterization} by the map to \((0, r) \times
    \sphere\) or \((0,r] \times \sphere\) and perform the same
    arguments for the extension to \(0\). Since \citet{guellaSchoenbergsTheoremPositive2019} consider
    \(\domain \times \sphere\) for general sets \(\domain\), their results still apply.
\end{remark}

 	\section{Strict positive definiteness}
\label{sec: strict positive definite}

Recall that a \shortkernel{} \(\kernel\) on $\real^\dims$ is called strictly positive definite, if for any $m\in\nat$, distinct $x_1,\dots ,x_m\in\real^\dims$ and $c\in \complex^m\setminus\{0\}$ we have 
\begin{align*}
	\sum_{i,j=1}^m c_i \conj{c_j}\kernel(x_i, x_j) > 0.
\end{align*}
In this section we characterize which continuous, isotropic
\pdkernels{} have this stronger property.
\begin{theorem}[Characterization of continous isotropic \spdkernels]
    \label{thm:strictly_pos_def}
    Let \(\dims \in \{3,\dots, \infty\}\).  For a continuous isotropic \pdkernel{} \(\kernel\colon \real^\dims
    \times \real^\dims \to \complex\)  with representation as in Theorem~\ref{thm: isotropy
    characterization}, the following assertions are \textbf{equivalent}:
    \begin{enumerate}[label={(\roman*)}]
        \item\label{it: strict positive definite}
        The \shortkernel{} \(\kernel\) is strictly positive definite.
        \item\label{it: infinitely even and odd elements in set}
        \(\alpha_0^{\dims}(0,0)>0\), and
        for any \(m\in \nat\), distinct \(r_1,\dots, r_m\in (0,\infty)\) and \(c \in \complex^m\setminus\{0\}\)
        the set
        \[
            \Bigl\{
                n \in \nat_0 :
                c^\transpose \Bigl[\alpha^{\dims}_n(r_i, r_j)\Bigr]_{i,j=1}^m \conj{c} > 0
            \Bigr\}
        \]
        contains infinitely many even and infinitely many odd integers.
        \item\label{it: even and odd kernels strict pos. def.}
        \(\alpha_0^{\dims}(0,0)>0\), and for each \(\gamma \ge 0\) the even and odd \shortkernels{}
        \[
            \alpha_\gamma^e(r,s) := \sum_{2n \ge \gamma} \alpha_{2n}^{\dims}(r,s)
            \quad\text{and}\quad
            \alpha_\gamma^o(r,s) := \sum_{2n+1 \ge \gamma} \alpha_{2n+1}^{\dims}(r,s)
        \]
        are both strictly positive definite on \((0,\infty)\).
    \end{enumerate}
    Furthermore, the following statement is \textbf{sufficient} for strict positive definiteness:
    \begin{enumerate}[label={(\roman*)},resume]
        \item\label{it: infinitely many even and odd n strict pos. def.}
        \(\alpha_0^{\dims}(0,0)>0\), and \(\alpha_n^{\dims}\) is strictly positive
        definite on \((0,\infty)\) for infinitely many even and odd \(n\in \nat_0\).
    \end{enumerate}
    If \(\alpha_0^{\dims}(0,0) =0\) but all other requirements in \ref{it:
    infinitely even and odd elements in set}, \ref{it: even and odd kernels
    strict pos. def.}, \ref{it: infinitely many even and odd n strict pos. def.}
    are satisfied, then \(\kernel\) is strictly positive definite on
    \(\real^\dims\setminus\{0\}\). \ref{it: infinitely even and odd elements in
    set} and \ref{it: even and odd kernels strict pos. def.} without \(\alpha_0^{\dims}(0,0) > 0\)
    are also necessary for strict positive definiteness on \(\real^\dims\setminus\{0\}\).
\end{theorem}

A proof can be found in Section \ref{sec:all-proofs}.
For the proof of Theorem \ref{thm:strictly_pos_def} we again utilize the map
from \((0,\infty)\times \sphere\) to \(\real^\dims\setminus\{0\}\) and base our
results on those by
\citet{guellaSchoenbergsTheoremPositive2019}.
Note that we consider \(\sphere\) while they considered \(\sphere[\dims]\),
leading to an index shift. Since they had to single out
the case \(\dims=2\), we do so as well and the following result covers this case. Observe that the first two conditions in
Theorem~\ref{thm:strictly_pos_def_dim2} are only necessary and not sufficient
for strict positive definiteness for \(\dims=2\) \citep[cf.][]{guellaPositiveDefiniteMatrixvalued2025}.

\begin{theorem}[Isotropic strictly positive definite kernels for \(\dims=2\)] 
    \label{thm:strictly_pos_def_dim2}
    Let \(\kernel\colon\real^2 \times\real^2 \to \complex\) be a continuous isotropic \pdkernel{} with decomposition as in Theorem~\ref{thm:
    isotropy characterization}. Then the following conditions are \textbf{necessary} for
    strict positive definiteness:
    \begin{enumerate}[label={(\roman*)}]
        \item \(\alpha_0^{2}(0,0)>0\), and
        for any \(m\in \nat\), distinct \(r_1,\dots, r_m\in (0,\infty)\) and \(c \in \complex^m\setminus\{0\}\)
        the set
        \[
            \Bigl\{
                n \in \integer :
                c^\transpose \Bigl[\alpha^{2}_{|n|}(r_i, r_j)\Bigr]_{i,j=1}^m \conj{c} > 0
            \Bigr\}
        \]
        intersects every full arithmetic progression\footnote{
            Recall, that a full arithmetic progression is a set of the form \((k \integer + m)\) with \(k\in \nat\) and \(m\in \{0,\dots k-1\}\).
        } in \(\integer\).

        \item \(\alpha_0^{2}(0,0)>0\), and for each full arithmetic progression \(S\) in \(\integer\)
        the \shortkernel{}
        \[
            (r,s) \mapsto \sum_{n\in S} \alpha^{2}_{|n|}(r,s) \qquad r,s\in (0,\infty).
        \]
        is strictly positive definite.
    \end{enumerate}
    A \textbf{sufficient} condition for strict positive definiteness of \(\kernel\)
    is
    \begin{enumerate}[label={(\roman*)},resume]
        \item \(\alpha_0^{2}(0,0)>0\), and the set \(\{n \in \integer : \alpha_{|n|}^{2} \text{ is strictly positive definite}\}\)
        intersects every full arithmetic progression in \(\integer\).
    \end{enumerate}
\end{theorem} 	\section{Examples of isotropic kernels}\label{sec: applications}

A survey of common kernel functions, whether 
implemented in scientific software libraries, e.g.\ \verb|KernelFunctions.jl|
\citep{galy-fajouJuliaGaussianProcessesKernelFunctionsjlV010652025},
or presented in foundational textbooks \citep[e.g.][Ch.\
4]{rasmussenGaussianProcessesMachine2006}, reveals that the vast
majority of canonical examples belongs to one of three categories
\begin{enumerate}[noitemsep]
	\item \textbf{stationary isotropic kernels}, i.e.\ functions
	of the distance \(\|x-y\|\),
	
	\item \textbf{dot product kernels}, i.e.\ functions of the inner product \(\langle x,y\rangle\),

	\item \textbf{specialized kernels} tailored for specific applications such as the
	\emph{neural network kernels}.
\end{enumerate}
Isotropic kernels jointly generalize the first two categories and
cover some specialized kernels, specifically the neural network kernels.
They can therefore serve as a unifying framework to avoid the separate treatment of
stationary isotropic kernels and dot product kernels
as in \citep{roosHighDimensionalGaussianProcess2021,amentScalableFirstOrderBayesian2022}.

\paragraph*{Outline of this section}
In Section \ref{sec: stationary isotropic and dot product} we illustrate the
unification of stationary isotropic and dot product kernels with examples.
In Section \ref{sec: machine learning} we highlight the importance of isotropic kernels
for the theory of machine learning. We demonstrate that limiting neural network kernels
are isotropic but not stationary and therefore do not fall in any of the previous
kernel classes while our results do apply. In Section \ref{sec: other examples}
we present two more examples of non-stationary isotropic kernels.

\subsection{Stationary isotropic and dot product kernels}
\label{sec: stationary isotropic and dot product}

Isotropic kernels include both stationary isotropic kernels and dot product
kernels. Hence both classes can be expressed by Theorem \ref{thm: isotropy characterization} as
\[
	\kernel(x,y)=\sum_{n=0}^\infty \alpha_n^\dims(\norm{x},\norm{y}) \uPnorm_n\bigl(\scp{\tfrac{x}{\norm x}}{\tfrac{y}{\norm y}}\bigr),
	\qquad \lambda = \tfrac{\dims-2}2,
\]
for specific choices of the kernels \(\alpha_n^\dims\). In this section
we derive these \(\alpha_n^\dims\) for the stationary isotropic and dot product kernels.

\begin{theorem}[Stationary isotropic kernels]
	\label{thm: stationary isotropic kernels}
	Recall that the stationary isotropic kernels in \(\real^\dims\) have been
	characterized by \citet{schoenbergMetricSpacesPositive1938} (see Table \ref{table: characterizations})
	to be of the form
	\[
		\kernel(x,y) = \begin{cases}
			\displaystyle
			\int_{[0, \infty)} \exp\bigl(-s^2\tfrac{\|x-y\|^2}2\bigr)d\mu(s)
			& \dims = \infty
			\\[3ex]
			\displaystyle
			\int_{[0, \infty)} \Omega_\dims(s\norm{x-y}) d\mu(s) & \dims < \infty,
		\end{cases}
		\qquad
	\]
	for a finite `Schoenberg measure' \(\mu\) on the interval \([0, \infty)\)
	and the characteristic function \(\Omega_\dims(t) = \E[e^{it\scp{v}{U}}]\)
	of a uniformly distributed \(U\sim\mathrm{Unif}(\sphere)\) for an arbitrary
	vector \(v\) with \(\norm{v}=1\). Then, the kernels \(\alpha_n\) in the
	isotropic kernel representation of Theorem \ref{thm: isotropy
	characterization} are given by
	\begin{equation}
		\label{eq: representations of alpha_n for stationary isotropic kernels}	
		\alpha_n^\dims(\norm{x}, \norm{y}) = \int_{[0, \infty)} \phi^\dims_n(s\norm{x})\phi^\dims_n(s\norm{y}) d\mu(s),
	\end{equation}
	with feature maps
	\[
		\phi^\dims_n(t) := \begin{cases}
			\displaystyle
			\exp(-\tfrac{t^2}2)\tfrac{t^n}{\sqrt{n!}} & \dims = \infty
			\\[1ex]
			\displaystyle
			2^\lambda \Gamma(\tfrac{\dims}{2}) \sqrt{N_n^\dims}\frac{\bessel_{n+\lambda}(t)}{t^\lambda} & \dims < \infty,
		\end{cases}
	\]
	where \(\bessel_\alpha\) is the Bessel function of the first kind \citep[9.1.10]{abramowitzHandbookMathematicalFunctions1964}
	and \(N_n^\dims = \dim \sphHarmonic_n\) is the dimension of the space of
	spherical harmonics \(\sphHarmonic_n\) of degree \(n\) in \(\sphere\) (cf.~Section \ref{sec: spherical harmonics}, Equation \eqref{eq: representations of spherical harmonics dimension}).
\end{theorem}
\begin{remark}
	Note that Squared Exponential, aka Gaussian, kernels correspond to Dirac Schoenberg measures \(\mu\).
\end{remark}
\begin{proof}[Proof of Theorem \ref{thm: stationary isotropic kernels}]
	Expanding the squared norm and using the series representation of the exponential results in
	\begin{align*}
		\kernel(x,y)
		&= \int_{[0, \infty)} \exp\bigl(-s^2\tfrac{\|x-y\|^2}2\bigr)d\mu(s)
		\\
		&= \int_{[0, \infty)}  \exp\bigl(-s^2\tfrac{\|x\|^2+\|y\|^2}2\bigr) \sum_{n=0}^\infty\frac{(s^2\langle x,y\rangle)^n}{n!} \mu(ds)
		\\
		&= \sum_{n=0}^\infty \underbrace{\Bigl(\int_{[0, \infty)} \exp\bigl(-s^2\tfrac{\|x\|^2+\|y\|^2}2\bigr) \tfrac{(s\|x\|)^n(s\|y\|)^n}{n!} \mu(ds)\Bigr)}_{= \alpha_n(\|x\|, \|y\|)}
		\bigl\langle \tfrac{x}{\|x\|},\tfrac{y}{\|y\|} \bigr\rangle^n.
	\end{align*}
	Note that the terms converge absolutely as can be seen by the application of Cauchy-Schwarz to \(\abs{\scp{x}y}^n\)
	in the second line and Hölder's inequality.
	In finite dimension we have by Lemma \ref{lem: characteristic function of uniform distribution on sphere}
	and \eqref{eq: standard bessel representation}
	\[
		\Omega_\dims(s\norm{x-y})
        = \sum_{n=0}^\infty N_n^\dims (2^\lambda \Gamma(\tfrac{\dims}2))^2 \frac{\bessel_{n+\lambda}(s\norm{x})}{(s\norm{x})^\lambda}\frac{\bessel_{n+\lambda}(s\norm{y})}{(s\norm{y})^\lambda} \uPnorm_n\bigl(\scp[\big]{\tfrac{x}{\norm{x}}}{\tfrac{y}{\norm{y}}}\bigr).
	\]
	The claim now follows directly from the representation of \(\kernel\).
\end{proof}

\begin{example}[Dot product kernels]\label{example:dot_product_kernels}
	Dot product kernels are of the form \(\kernel(x,y):= \kappa(\scp{x}{y})\) and
	hence are isotropic. Restricted to the sphere they exhaust all isotropic
	kernels and are characterized by \citet{schoenbergPositiveDefiniteFunctions1942} (see also Table \ref{table: characterizations}).
	A dot product kernel $\kernel$ on the $\ell_2$-sphere
is of the form
	\begin{equation}
		\label{eq: representation dot product kernels}	
		\kernel(x,y) = \sum_{n=0}^\infty a_n \langle x,y\rangle^n,
	\end{equation}
	with \((a_n)_{n\in \nat_0}\) a nonnegative, absolutely summable sequence.  
	Such kernels are also
	referred to as `spherical spin glasses' in statistical
	physics, e.g.\@ \citep[][]{panchenkoSherringtonKirkpatrickModel2013}.

	The dot product kernels on the entire sequence space \(\ell^2\) are
	still characterized by \eqref{eq: representation dot product kernels}, together with the simple requirement that \((a_nt^n)_{n\in \nat_0}\) must be summable for any \(t>0\), see
	\citep{pinkusStrictlyPositiveDefinite2004} and references therein.
	A simple alternative proof of this result based on the characterization of
	isotropic kernels (Theorem \ref{thm: isotropy characterization}) can be
	found in the appendix (Lemma \ref{lem: classification of dot product kernels
	infinite dim}).
\end{example}

\subsection{Machine Learning}\label{sec: machine learning}

Non-stationary isotropic kernels arise naturally in the analysis of infinitely
wide neural networks. These limiting networks offer tractable surrogates for
highly overparametrized neural networks. Specifically, under suitable random
initialization, the infinite width limit of a neural network is a Gaussian
process,
see e.g.\
\citep[][]{williamsComputingInfiniteNetworks1996,choKernelMethodsDeep2009,leeDeepNeuralNetworks2018,yangTensorProgramsWide2019},
and asymptotic training dynamics can be described by kernel gradient descent governed by the \emph{Neural Tangent Kernel (NTK)} \citep[][]{jacotNeuralTangentKernel2018,yangTensorProgramsII2020}.
\citet[Section F]{benningRandomFunctionDescent2024} demonstrate using
a simple linear toy model that
stationarity in contrast to isotropy is an unrealistic assumption for the
risk functions encountered in machine learning. 
In this section we explain that both the covariance kernel of the
infinite-width Gaussian process as well as the NTK of feed-forward networks are non-stationary isotropic
kernels, and the theory developed here applies to them. Furthermore, we discuss in Remark \ref{rem:other-architectures} that this is typical for many architectures beyond feed-forward networks. The class of non-stationary isotropic kernels is therefore
a natural and central class of kernels in machine learning.

\begin{definition}[Feed-forward neural network]\label{def:feed-forward-net} An \(L\)-layer neural network
\(f_\theta\colon \real^\dims \to \real\)
is defined as \(f_\theta(x) := z^{(L)}\), where \(x^{(0)} := x\) and the
	collection of all parameters \(\theta=(W^{(l)}, b^{(l)})_{l=1,\dots, L}\)
	recursively gives
	\begin{alignat*}{2}
z^{(l)} & := \gamma_lW^{(l)} x^{(l-1)} +  b^{(l)} &&\in \real^{n_l},
		\\ x^{(l)} &:= \phi(z^{(l)}) &&\in \real^{n_l}.
\end{alignat*}
	The \emph{activation function} $\phi\colon\real\to\real$ is applied
	coordinate-wise and the normalization factors
	\[
		\gamma_l := \begin{cases}
			\frac{1}{\sqrt{n_{l-1}}} & l \ge 2
			\\
			1 & l=1
		\end{cases}
	\]
	allow for a non-trivial limit in the width of hidden layers.
	
	In the \emph{NTK initialization,} all parameters \(\theta\) are independently standard normal distributed.
\end{definition}

The limit kernels of a feed-forward network naturally depend on the activation function. However in any case they are isotropic.
\begin{corollary}[Convergence to isotropic Gaussian processes]
\label{cor: limit_kernels}
	Consider a feed-forward network \(f_\theta\) with NTK initialization of the
	parameters \(\theta\) as in Definition \ref{def:feed-forward-net} with polynomially
	bounded activation function $\phi$. Assume for simplicity that all hidden
	layers have the same width
	\(n=n_l\). Then there exists a Gaussian process $g$ on $\real^\dims$ such that
	as the widths $n$ of the hidden layers increase, the neural
	networks converge in finite-dimensional distributions to $g$, that is
	\[
	f_\theta \overset{\mathrm{f.d.d.}}{\longrightarrow} g
	\qquad (n\to \infty).
	\]
	This process $g$, the so called \emph{Neural Network
		Gaussian process (NNGP)}, has mean zero and covariance \shortkernel{}
		$\nngp{}$, which is an isotropic kernel (typically non-stationary)
		of the form
	\begin{equation}
		\label{eq: joint nngp form}	
		\nngp{} (\xxs,\xys)
		= \sum_{m=0}^\infty \alpha_m(\norm {\xxs}, \norm{\xys})
		\bigl\langle \tfrac{\xxs}{\|\xxs\|}, \tfrac{\xys}{\|\xys\|}\bigr\rangle^m.
	\end{equation}
	Here \(\alpha_m\colon [0,\infty)^2\to\real\) are unique, continuous,
	\pdkernels. These \shortkernels{} are summable, i.e.
	\(\sum_{m=0}^\infty \alpha_m(r,r) < \infty\) for any \(r\in
	[0,\infty)\), and zero at the origin, i.e. 
	\(
	\alpha_m(r,0) = \alpha_m(0,r) = 0
	\)	
	for all $ r \in [0,\infty), m\ge 1$.
\end{corollary}
\begin{proof}
	Following \citet{haninRandomNeuralNetworks2023}, we can combine Eq.\ (1.7) and Eq.\ (1.8) therein to obtain
	\begin{align*}
		\nngp{l}(\xxs,\xys)
		&= 1+ F_\phi\Bigl(
		\nngp{l-1}(\xxs,\xys),\nngp{l-1}(\xxs,\xxs),\nngp{l-1}(\xys,\xys)
		\Bigr)
		\quad (l\ge 2)\\
		\nngp{1}(\xxs,\xys) &= 1  + \scp {\xxs}{\xys}.
	\end{align*}
	The exact form of $F_\phi$ is irrelevant here; it
	suffices to note that it is continuous, which is a consequence of the integral representation in \citep[Eq.\ (1.8)]{haninRandomNeuralNetworks2023}.  
	This recursion demonstrates the isotropy of $\nngp{} =\nngp{L}$ as
	it only depends on $\scp{\xxs}{\xys}$, $\norm{\xxs}^2$ and $\norm{\xys}^2$ 
	which are invariant under mutual linear isometries. 
Moreover, the function $F_\phi$ does not depend on the input dimension and hence the kernel function
		\(\kappa\) with
		\[
		\kappa(\|x\|, \|\xys\|, \langle x,\xys\rangle) \equalDef \nngp{}(x,\xys)
		\]	
	is valid in all finite dimensions. Hence, it produces a positive definite kernel on \(\ell^2\) by Lemma \ref{lem: valid in all dimensions}.
	The representation \eqref{eq: joint nngp form} of \(\nngp{}\) follows from Theorem \ref{thm: isotropy characterization}.
\end{proof}

\begin{remark}[Neural Tangent Kernel]
	From similar recursive representations it is possible to infer isotropy
	for the \emph{Neural Tangent Kernel} (NTK) \cite{jacotNeuralTangentKernel2018}
	defined by
	\[
		\ntk{}(\xxs,\xys):= \lim_{n\to \infty} \scp{\nabla_\theta f_\theta (\xxs)}{\nabla_\theta f_\theta(\xys)},
	\]
	which describes training dynamics of highly overparameterized networks.
\end{remark}
The isotropy of the limit neural network \shortkernels{} stems from the
isotropic processing of the inputs in the first hidden layer by the linear
multiplication with the Gaussian matrix $W^{(1)}$.
Specifically, let
${W_{i,\cdot}^{(1)}} \sim \mathcal{N}(0,\identity_\dims)$, then  for inputs
$\xx, \xy \in\real^\dims$ the
post-activation vectors have the coordinate-wise distribution
\begin{align}
	\label{eq:linear_layer_isotropy}
	\begin{pmatrix}
		W^{(1)}_{i,\cdot}\xx \\
		W^{(1)}_{i,\cdot} \xy
	\end{pmatrix} 
	\sim \mathcal{N}\Big(0,\begin{pmatrix}
		\scp {\xx}{\xx} & \scp{\xx}{\xy}\\
		\scp{\xy}{\xx} & \scp{\xy}{\xy}
	\end{pmatrix}\Big)
\end{align}
and the coordinates are independent.
This is an isotropic distribution in the sense that it is invariant under $A\in O(\dims)$ acting simultaneously on $\xx$ and $\xy$.
The isotropy \eqref{eq:linear_layer_isotropy} of linear layers, which are
central in neural network architectures, leads to isotropic limit
\shortkernels{} in much more general settings:

\begin{remark}[Other Network Architectures]\label{rem:other-architectures}
	Using `Tensor programs' \citep{yangTensorProgramsII2020} it is possible to
	recursively calculate the limiting kernels of numerous neural network
	architectures such as
	RNNs 
and 
	transformers as well as networks containing batchnorm and
	pooling layers. 
	Remarkably, it turns out that any network where the input is first processed by a fully connected linear layer (as in Def.\@ \ref{def:feed-forward-net}) yields isotropic limit kernels.
	The recursive formulas obtained tend to be
	complicated, but it is often straightforward to infer isotropy. 
\end{remark}

Neural network kernels can appear to be entirely unrelated, and the following examples illustrate the variety of such kernels.
What connects them
is that they are isotropic and thereby of the form \eqref{eq: joint nngp form}. 
\begin{example}[NNGP-kernels of shallow networks]	
	\leavevmode
	Consider a neural network as in Definition
	\ref{def:feed-forward-net} with $L=2$, but \emph{without bias }$b^{(l)}$.
	Varying the activation function $\phi$ yields a wide range of NNGP-kernels.
	Note that the following formulas differ slightly from the cited works since our
	normalization $\gamma_l$ only coincides with their normalization up to constants.
	\begin{itemize}		
		\item For the ReLU activation function \(\phi(x) = \max\{0,x\}\)		
		\citet{choKernelMethodsDeep2009} calculate the NNGP-kernel to be the `arccosine kernel'		
		\[
			\nngp{}(\xxs,\xys) = 
			\frac1{2\pi}\|x\| \|\xys\| J(\cos^{-1}(\scp[\big]{\tfrac{x}{\norm{x}}}{\tfrac{\xys}{\norm{\xys}}})),\\			
		\]
		with \(J(\theta) \coloneq \sin(\theta) + (\pi -\theta)\cos(\theta)\).
		\item 	\citet{williamsComputingInfiniteNetworks1996} shows that 		
		the NNGP-kernel for a shallow neural network with the `error function' $\phi(x)=\tfrac{2}{\sqrt \pi}\int_0^t e^{-s^2}d s$ as		
		activation function is the `arcsin kernel'		
		\[
		\nngp{} (\xxs,\xys) = 
		\frac {2}{\pi}\arcsin\Bigl(\frac{2\scp x{\xys}}{\sqrt{(1+2\|x\|^2)(1+2\|\xys\|^2)}}\Bigr).
		\qedhere
		\]
	\end{itemize}
\end{example}

\paragraph*{Characterization of shallow NNGP-kernels on the sphere}
For the class of activation functions
\[
	\phi(x) = \sum_{n=0}^\infty a_n h_n(x),
\]
where \(a_n \in \real\) and $h_n$ denotes the $n$-th Hermite polynomial, \citet[Thm.\@
3.1]{simonReverseEngineeringNeural2022} show that the NNGP kernel of a \(2\)-layer neural network for spherical inputs, i.e.\@ $\norm x=\norm {\xys}=1$,
is of the form
\[
	\nngp{}(x,y) = \sum_{n=0}^\infty a_n^2 \scp x{\xys}^n.
\]
Schoenberg's characterization of isotropic kernels on the sphere
\citep{schoenbergPositiveDefiniteFunctions1942} (Table \ref{table: characterizations}) then implies that \emph{any}
continuous isotropic kernel on $\sphere[\infty]$ may be interpreted as an
NNGP-kernel with spherical input.

For future research, it is interesting to investigate whether every non-stationary isotropic
kernel on \(\ell^2\) can be realized as an NNGP-kernel outside of the sphere.
Our characterization of isotropic kernels is a first step towards answering this question.

\subsection{Other examples of non-stationary isotropic kernels}
\label{sec: other examples}

\begin{example}[Gegenbauer generating kernel]
	Using the generating function of the Gegenbauer polynomials \citep[(2.3)]{reimerGegenbauerPolynomials2003}
	\[
		\sum_{n=0}^\infty z^n\uP_n(t)  = (1-2tz+z^2)^{-\lambda}, \qquad |z|<1,\, t\in [-1,1],\, \lambda \neq 0,
	\]
	any continuous kernel \(\beta \colon [0, \infty)^2 \to \real\) with \(\beta(r, 0) = 0\) and \(\abs{\beta(r,r)} < 1\) for all \(r\in [0,\infty)\)
	naturally yields an isotropic kernel on \(\real^\dims\) for \(\dims >2\) given by
	\begin{align*}
		\kernel(x,y)
		&\coloneq \Bigl(1-2\scp[\big]{\tfrac{x}{\norm{x}}}{\tfrac{y}{\norm{y}}}\beta(\norm{x},\norm{y})+\beta(\norm{x},\norm{y})^2\Bigr)^{-\lambda}
		&\lambda = \tfrac{\dims-2}{2}
		\\
		&= \sum_{n=0}^\infty \beta(\norm{x},\norm{y})^n \uP_n\bigl(\scp[\big]{\tfrac{x}{\norm{x}}}{\tfrac{y}{\norm{y}}}\bigr)
		&x,y \in \real^\dims.
	\end{align*}
	In the notation of Theorem \ref{thm: isotropy characterization}, the decomposing kernels are therefore essentially
	powers of \(\beta\) up to normalization, that is \(\alpha_n^\dims(r,s) = \beta(r,s)^n \uP_n(1)\).
\end{example}

Although kernels defined through hyperbolic distance are not translation
invariant in the Euclidean coordinates of the Poincaré ball, they may still be
isotropic with respect to Euclidean rotations. The following example shows how
our isotropy characterization can be used to prove positive definiteness of such
a hyperbolic distance kernel.

\begin{example}[Hyperbolic distance kernel]
	Let \(\metric\) denote the hyperbolic distance on the Poincaré ball, that is
	\[
		\metric(x,y) = \arcosh\Bigl(1+2\frac{\norm{x-y}^2}{(1-\norm{x}^2)(1-\norm{y}^2)}\Bigr), \qquad \norm{x},\norm{y} < 1.
	\]
	The Gaussian kernel \(\exp(-\gamma \metric(x,y)^2)\) is never positive definite in hyperbolic space \citep{dacostaInvariantKernelsRiemannian2025},
	and \citet{dacostaInvariantKernelsRiemannian2025} argue that a suitable replacement for the Gaussian kernel in hyperbolic space is given by the kernel
	\[
		\kernel_\gamma(x,y)
		= \cosh(\metric(x,y))^{-\gamma} 
		\qquad \text{for}\qquad \gamma>0.
	\]
	This kernel has the following alternative representations
	\[
		\kernel_\gamma(x,y)
		= \Bigl[\frac{(1-\norm{x}^2)(1-\norm{y}^2)}{(1+\norm{x}^2)(1+\norm{y}^2) - 4 \scp{x}{y}}\Bigr]^\gamma
		= \sum_{n=0}^\infty \alpha_n(\norm{x}, \norm{y}) \scp[\big]{\tfrac{x}{\norm{x}}}{\tfrac{y}{\norm{y}}}^n,
	\]
	with \(
		\alpha_n(\norm{x}, \norm{y})
		\coloneq \phi_n(\norm{x})\phi_n(\norm{y})
	\)
	defined by the real-valued feature maps
	\[
		\phi_n(r) \coloneq \sqrt{\tfrac{(\gamma)_n}{n!}}\bigl(\tfrac{1-r^2}{1+r^2}\bigr)^\gamma \bigl(\tfrac{2r}{1+r^2}\bigr)^n, \qquad r\in [0,1).
	\]
	Here \((\gamma)_n = \gamma(\gamma+1)\cdots(\gamma+n-1)\) is the Pochhammer
	symbol.  Remark \ref{rem: isotropic kernels on balls} then implies that
	\(\kernel_\gamma\) is a positive definite kernel on the unit ball in
	\(\real^\dims\) for any \(\dims\ge 2\).
	
\end{example}

\begin{proof}
	The first representation follows directly from the definition of the
	hyperbolic distance and the expansion \(\norm{x-y}^2 = \norm{x}^2 +
	\norm{y}^2 - 2\scp{x}{y}\). For the second representation write \(r\coloneq \norm{x}\), \(s\coloneq \norm{y}\)
	and \(t\coloneq \scp{\frac{x}{\norm{x}}}{\frac y{\norm{y}}}\). Then
	\[
		\kernel_\gamma(x,y)
		= \Bigl[\frac{(1+r^2)(1+s^2) - 4 r s t}{(1-r^2)(1-s^2)}\Bigr]^{-\gamma}
		= b_\gamma(r)b_\gamma(s)\bigl[1-q(r)q(s)t\bigr]^{-\gamma}
	\]
	with \(q(r)\coloneq \tfrac{2r}{1+r^2}\) and \(b_\gamma(r) \coloneq \bigl(\tfrac{1-r^2}{1+r^2}\bigr)^\gamma\).
	Rewriting the binomial series \citep[3.6.8]{abramowitzHandbookMathematicalFunctions1964}
	\[
		(1-x)^{-\gamma} = \sum_{n=0}^\infty \binom{-\gamma}{n}(-x)^n
	\]
	with
	\(
		\binom{-\gamma}{n}
		= (-1)^n \frac{\gamma(\gamma+1)\cdots(\gamma+n-1)}{n!}
		= (-1)^n \frac{(\gamma)_n}{n!}
	\)
	yields
	\[
		[1-q(r)q(s)t]^{-\gamma} = \sum_{n=0}^\infty \frac{(\gamma)_n}{n!}q(r)^n q(s)^n t^n
	\]
	and therefore
	\[
		\kernel_\gamma(x,y)
		= \sum_{n=0}^\infty \phi_n(r)\phi_n(s) t^n
		\quad\text{with}\quad
		\phi_n(r) \coloneq \sqrt{\tfrac{(\gamma)_n}{n!}} b_\gamma(r) q(r)^n.
	\]
	Resubstituting the definitions of \(q\) and \(b_\gamma\) yields the claimed representation. 
\end{proof}

 	\section{Application: Isotropic Gaussian random function representation}
\label{sec: GRF representation}
For a positive definite kernel \(\kernel\) there always exists a Gaussian random
function (GRF) \(\rf\), also referred to as a Gaussian process or Gaussian
random field, with covariance
\[
    \Cov(\rf(x), \rf(y)) = \kernel(x,y)
\]
by Kolmogorov's extension
theorem \citep[e.g.][Thm.~14.36]{klenkeProbabilityTheoryComprehensive2014}. 
Naturally, such an existence result is highly abstract and more explicit
representations are desirable in both theory and practice. 
The mere existence of the GRF does not suggest a method for simulation for example.
Using the previously established characterization of isotropic positive definite kernels $\kernel$ in Theorem \ref{thm: isotropy characterization}, we obtain an explicit representation of the corresponding GRFs.

\begin{theorem}[Representation of isotropic Gaussian random functions]
    \label{thm: representation isotropic random functions}
    Let
    \[
        \kernel(x,y) = \sum_{n=0}^\infty \alpha_n^\dims \bigl(\|x\|, \|y\|\bigr) \uPnorm_n\bigl(\bigl\langle \tfrac{x}{\|x\|}, \tfrac{y}{\|y\|}\bigr\rangle\bigr),
        \qquad x,y\in \real^\dims,\; \lambda = \tfrac{\dims-2}{2},
    \]
    be a real-valued, isotropic, continuous \pdkernel{}, with the unique
    one dimensional kernels \(\alpha_n^\dims\) from the characterization of
    Theorem \ref{thm: isotropy characterization}. A centered Gaussian random
    function \(\rf\) with continuous covariance $\kernel$
    may be represented as
    \[
        \rf(x) \coloneq \sum_{n=0}^\infty \rf_n(x),
    \]
    using independent Gaussian random functions \(\rf_n\) explicitly defined below,
    where the series converges for every \(x\in \real^\dims\) both in \(L^2(\Omega)\) and almost surely.
    \begin{enumerate}
        \item \textbf{Case \( \dims=\infty\):}
        Let \(\rg_0\) be a real-valued, centered Gaussian random function on \([0, \infty)\) with covariance kernel \(\alpha_0^\infty\)
        and \(\rg_n^{i_1, \dots, i_n}\) with \(i_j \in \nat\) real-valued, centered, independent Gaussian random functions
        on \([0, \infty)\) with  covariance kernel \(\alpha_n^\infty\).
        Then we define
        \begin{equation}
            \label{eq: representation of GRF valid in all dimensions}            
            \rf_n(x) := \!\!\sum_{i_1,\dots, i_n=1}^\infty \!\! \rg_n^{i_1, \dots, i_n}(\|x\|) \hat x_{i_1}\cdot \hdots \cdot \hat x_{i_n},
            \qquad 
            \rf_0(x) := \rg_0(\|x\|),
        \end{equation}
        using the notation \(\hat x \coloneq x/\|x\|\) for \(x\neq 0\) and \(\hat 0 \coloneq 0\).
        The series in \eqref{eq: representation of GRF valid in all dimensions} converges in \(L^2(\Omega)\) and almost surely for each \(x\in \ell^2\).

        \item \textbf{Case \(\dims<\infty\):} Let \(Y_{n,k}\) be an orthonormal basis of the
        space of spherical harmonics \(\sphHarmonic_n\) of degree \(n\in \nat\)
        in \(\sphere\), and for \(k\in\set{1,\dots, N_n^\dims}\) with
        \(N_n^\dims = \dim \sphHarmonic_n\) let \(\rg_n^k\) be real-valued,
        centered independent Gaussian random functions on \([0, \infty)\) with
        covariance kernel \(\alpha_n^\dims\).
        Then we define \(\rf_0(x) \coloneq \rg_0^1(\|x\|)\) and for \(n\ge 1\), \(\rf_n(0) \coloneq 0\)  and
        \begin{equation}
            \label{eq: spherical harmonics representation of GRF}    
            \rf_n(x) \coloneq 
                \frac1{\sqrt{N_n^\dims}}\sum_{k=1}^{N_n^\dims} \rg_n^k(\|x\|) Y_{n,k}\bigl(\tfrac x{\norm{x}}\bigr)
            \qquad x\in \real^\dims\setminus\{0\}.
        \end{equation}
    \end{enumerate}
\end{theorem}
Random functions whose distribution is invariant under rotation, that is \(\Pr_\rf =
\Pr_{\rf \circ U}\) for \(U\in O(\dims)\) are called isotropic.  The representation
above, characterizes all mean-squared continuous, centered, isotropic Gaussian random
functions, since Gaussian random functions are characterized by mean and
covariance and Theorem \ref{thm: isotropy characterization} characterizes all
continuous isotropic covariance kernels.

The infinite dimensional representation \eqref{eq: representation of GRF valid in all dimensions} generalizes the well-known series representation
of isotropic GRFs on the sphere in terms of polynomials of order \(n\) with random
coefficients.
These random polynomials of order \(n\) are referred to as ``pure
\(n\)-spins'' in the spin-glass
community \citep[e.g.][]{auffingerOptimizationRandomHighDimensional2023,rosHighdLandscapesParadigm2023,panchenkoSherringtonKirkpatrickModel2013}.
On the sphere, the finite dimensional representation is also well known \citep[e.g.][Thm.~5.13]{marinucciRandomFieldsSphere2011}.

\begin{remark}[Existence and continuity of \(\rg_n\)]
    The existence of the Gaussian random functions \(\rg_n\) in \eqref{eq: representation of GRF valid in all dimensions} and \eqref{eq: spherical harmonics representation of GRF} follows from Kolmogorov's extension theorem \citep[e.g.][Thm.~14.36]{klenkeProbabilityTheoryComprehensive2014} since
    the \(\alpha_n^\dims\) are positive definite kernels on \([0,\infty)\).
    By Theorem \ref{thm: isotropy characterization}, these kernels \(\alpha_n^\dims\) are furthermore continuous and the resulting process \(\rg_n\)
    therefore mean-square continuous \citep[Thm.~5.3.3]{scheurerComparisonModelsMethods2009}. 
    Sample path continuity requires slightly stronger assumptions such as Kolmogorov-Chentsov's criterion \citep[Thm.~21.6]{klenkeProbabilityTheoryComprehensive2014}
    or other criteria \citep[e.g.][Thm.~1.4.1]{talagrandRegularityGaussianProcesses1987,adlerRandomFieldsGeometry2007}.
\end{remark}

\begin{remark}[Simulation]
    A possible application of Theorem
    \ref{thm: representation isotropic random functions} is the simulation of
    isotropic random functions. 
    If the kernels \(\alpha_n^\dims\) allow a simple representation of the form \(\alpha_n^\dims(r,s) = \sum_{i=1}^m\phi_i(r)\phi_i(s)\) for some feature maps \(\phi_i\), then we can sample the GRF $\rg_n^k$ as
    \[
        \rg_n^k(r) \sim \sum_{i=1}^m \xi_{n,i}^k \phi_i(r) 
        \quad \text{with} \quad
        \xi_{n,i}^k\overset{\iid}\sim \normal(0,1),
    \]
    which allows approximating sampling of the GRF $\rf$ by truncating the series in \eqref{eq: representation of GRF valid in all dimensions}.
    This is the case in the following examples:
    \begin{enumerate}
        \item \textbf{Stationary isotropic kernels:}
        If the Schoenberg measure \(\mu\) (cf.~Theorem~\ref{thm: stationary isotropic kernels}) is supported on finitely many points \(\gamma_1,\dots, \gamma_m\),
        then by \eqref{eq: representations of alpha_n for stationary isotropic
        kernels}
        \[
            \alpha_n^\dims(r,s)
            = \sum_{i=1}^m \varphi_i(r)\varphi_i(s)
            \qquad \text{with} \quad
            \varphi_i(r) \coloneq \sqrt{\mu(\gamma_i)}\phi_n^\dims(\gamma_i r),
        \]
        where the \(\phi_n^\dims\) are defined below \eqref{eq: representations
        of alpha_n for stationary isotropic kernels}.

        \item \textbf{Dot product kernels:} Here \(\alpha_n^\dims(r,s) = r^n s^n\)
        and the one dimensional feature map is trivially given by \(\phi(r) = r^n\).
    \end{enumerate}
    \emph{Computational complexity.} For the simulation of a GRF \(\rf\) with covariance kernel
    \(\kernel\) in dimension $\dims\in \nat$ one may in principle use a higher dimensional representation in dimension ${\dims'}\in \nat\cup\{\infty\}$ using the canonical embedding $\real^\dims \subset \real^{\dims'}$,
    which appends zeros to a given vector $x\in \real^\dims$.
    Then there are only finitely many non-zero summands in \eqref{eq: representation of GRF valid in all dimensions} 
    and the infinite dimensional representation becomes a finite sum.
    However, the computational complexity of \eqref{eq: representation of GRF valid in all dimensions}
    is then \(\dims^n\) to simulate each \(\rf_n\), which is generally much more expensive than the
    complexity of \eqref{eq: spherical harmonics representation of GRF}, which
    is of order \(N_n^\dims \sim \tfrac{2}{(\dims-2)!}n^{\dims-2}\) as \(n\to
    \infty\) (cf.~Lem.~\ref{lem: spherical harmonics dimension}). Using the
    appropriate dimensional Gegenbauer representation is therefore more
    efficient computationally, assuming this representation is available.
\end{remark}

\begin{proof}[Proof of Theorem \ref{thm: representation isotropic random functions}]
    We start by calculating the covariance of the \(\rf_n\) in both cases.
    In the case \(n=0\) we have
    \[
        \Cov(\rf_0(x), \rf_0(y)) = \Cov(\rg_0(\|x\|), \rg_0(\|y\|)) = \alpha_0^\dims(\|x\|, \|y\|), 
    \]
    and if either \(x=0\) or \(y=0\), then the covariance is zero for \(n \ge 1\) as
    we have \(\rf_n(0)=0\) in this case. So we will assume \(x,y\neq 0\) and \(n \ge 1\) without
    loss of generality in the following.
    Starting with the case \(\dims < \infty\) we can eliminate the mixed terms
    using the independence of the centered \(\rg_n^k\), which yields
    \begin{align*}
        \Cov(\rf_n(x), \rf_n(y))
        &= \E[\rf_n(x) \rf_n(y)]
        \\
        &= \frac1{N_n^\dims}\sum_{k,l=1}^{N_n^\dims} \underbrace{\E[\rg_n^k(\|x\|) \rg_n^l(\|y\|)]}_{
            = \alpha_n^\dims(\|x\|, \|y\|) \delta_{kl}
        } Y_{n,k}(\hat x) Y_{n,l}(\hat y)
        \\
        &= \alpha_n^\dims(\|x\|, \|y\|) \underbrace{\frac1{N_n^\dims}\sum_{k=1}^{N_n^\dims}  Y_{n,k}(\hat x) Y_{n,k}(\hat y)}_{= \uPnorm_n(\scp{\hat x}{\hat y})}.
    \end{align*}
    The last equality follows from the addition formula for spherical harmonics (Theorem \ref{thm: addition theorem}).
    Heuristically, the case \(\dims=\infty\) can be proven the same way. But as we are
    dealing with an infinite series a formal argument requires more care.
    Observe that every summand \(X_{\mathbf i}(x) \coloneq \rg_n^{\mathbf i}(\norm{x})\prod_{i\in \mathbf{i}} \hat x_{i} \)
    with \(\mathbf i \coloneq (i_1,\dots, i_n){\in \nat^n}\) is a centered Gaussian random
    variable with variance
    \begin{equation}
        \label{eq: one point variance}
        \E\bigl[X_{\mathbf i}(x)^2\bigr]
        = \E\bigl[\rg_n^{\mathbf i}(\norm{x})^2\bigr] \prod_{i\in \mathbf{i}} \hat x_{i}^2
        = \alpha_n^\infty(\|x\|, \|x\|) \prod_{i\in \mathbf{i}} \hat x_{i}^2.
    \end{equation}
    For \(x,y\in \ell^2\) we have
    \begin{equation}
        \label{eq: sum over prod i x_i y_i}
        \sum_{\mathbf i \in \nat^n} \prod_{i\in \mathbf{i}} \hat x_{i} \hat y_{i}
        = \sum_{i_1,\dots, i_n=1}^\infty \!\!\hat x_{i_1}\cdot \hdots \cdot \hat x_{i_n} \hat y_{i_1}\cdot \hdots \cdot\hat y_{i_n}
        = \Bigl(\sum_{i=1}^\infty \hat x_i \hat y_i\Bigr)^n
        = \scp{\hat x}{\hat y}^n.
    \end{equation}
    In particular the variance in \eqref{eq: one point variance}
    is summable. And since \(\norm{\hat x} = 1\) we have
    \begin{equation}
        \label{eq: variance is summable}    
        \sum_{\mathbf i \in \nat^n} \E\bigl[X_{\mathbf i}(x)^2\bigr]
        = \alpha_n^\infty(\|x\|, \|x\|) \sum_{\mathbf i \in \nat^n} \prod_{i\in \mathbf{i}} \hat x_{i}^2
        = \alpha_n^\infty(\|x\|, \|x\|) < \infty.
    \end{equation}
    This implies the summability assumption of Lemma \ref{lem: series
    convergence} is satisfied and the infinite sum over
    \(X_{\mathbf i}(x)\) converges in \(L^2(\Omega)\) and almost surely to a
    centered Gaussian random variable \(\rf_n(x)\) with variance \(\alpha_n^\infty(\|x\|,
    \|x\|)\). Applying this argument for every \(x\in \ell^2\) yields a collection of
    Gaussian random variables \((\rf_n(x))_{x\in \ell^2}\). To prove that this collection is a Gaussian
    random function, we need to show that every finite dimensional marginal is
    Gaussian with the correct covariance structure. But this follows by
    considering the vectors \(X_{\mathbf i}=(X_{\mathbf i}(x^1), \dots,
    X_{\mathbf i}(x^m))\) for \(m\in \nat\), \(x^1,\dots, x^m\in \ell^2\) and all \(\mathbf i\in\nat^n\), which
    have the covariance matrix
    \[
        \Cov(X_{\mathbf i})_{kl}
        = \E[X_{\mathbf i}(x^k) X_{\mathbf i}(x^l)]
        = \alpha_n^\infty(\|x^k\|, \|x^l\|) \prod_{i\in \mathbf{i}} \hat x_{i}^k \hat x_{i}^l.
    \]
    By \eqref{eq: sum over prod i x_i y_i} these sum to
    \begin{equation}
        \label{eq: sum of covariance matrices}
        \sum_{\mathbf i \in \nat^n} \Cov(X_{\mathbf i})_{kl}
        = \alpha_n^\infty(\|x^k\|, \|x^l\|)\scp{\hat x^k}{\hat x^l}^n
\end{equation}
    Lemma \ref{lem: series convergence} implies that the covariance matrix of the series of \(X_{\mathbf i}\)
    is the series of covariance matrices, which thus yields that \(\rf_n\) is a
    centered Gaussian random function with covariance kernel in \(x^k\) and
    \(x^l\) given in \eqref{eq: sum of covariance matrices}.  Collecting all the
    cases, we have both in the finite and infinite dimensional case for all \(n\in \nat_0\)
    \[
        \Cov(\rf_n(x), \rf_n(y))
        = \alpha_n^\dims(\|x\|, \|y\|) \uPnorm_n(\scp{\hat x}{\hat y}) \eqcolon \kernel_n(x,y),
    \]
    where we recall \(\uPnorm_n(t) = t^n\) for \(\dims=\infty\) and \(\uPnorm_0 \equiv 1\). The \(\rf_n\) are thus
    independent centered Gaussian random functions with covariance kernels \(\kernel_n\).
    Since the infinite sum over \(\kernel_n\) converges absolutely to \(\kernel\) by
    Cauchy-Schwarz, a final application of Lemma \ref{lem: series convergence} yields the claimed convergence of
    the infinite sum over \(\rf_n\).
\end{proof} 
	\subsection*{Acknowledgements}
This work was partially supported by the Deutsche Forschungsgemeinschaft (DFG, German Research Foundation)  through the project STE 1074/5-1, within the DFG priority programme SPP 2298 “Theoretical Foundations of Deep Learning”, to M.D.\ Sch\"olpple. 	\bibliography{corrected_references}
	
	\appendix

 	\section{Proofs}
\label{sec: proofs}

In this section we provide the proofs omitted in the main text.
We begin with a recapitulation of well known results and tools on which our
theorems are built.

\subsection{Recapitulation of foundations}
\label{sec: recap foundations}
Our work is largely built on \citet{guellaSchoenbergsTheoremPositive2019}. 
Note that in contrast to them we work on $\sphere $ instead of $\sphere[\dims]$ and use \emph{normalized} Gegenbauer polynomials in contrast to unnormalized ones.
For the reader's convenience we therefore restate their main theorem in our notation and denote by $\domain$ some arbitrary non-empty set.

\begin{theorem}[Characterization of kernels on \(\domain \times \sphere\), \cite{guellaSchoenbergsTheoremPositive2019}]
	\label{thm: guella XxS characterization}
	For the dimension \(\dims \in \{2,3,\dots, \infty\}\) define \(\lambda:=\frac{\dims-2}2\). For a
	function \(\kernel: (\domain \times \sphere)^2 \to \complex\) which is
	isotropic on the sphere, i.e.\  of
	the form 
	\[
	\kernel((r,v), (s,w)) = f_{\sphere[]}(r,s, \langle v ,w\rangle),
	\]
	the following are equivalent 
	\begin{enumerate}[label=(\roman*)]
		\item The function \(\kernel \) is a \pdkernel{} and for all $r,s\in \domain$ the function
		\(\rho \mapsto f_{\sphere[]}(r,s, \rho)\) is continuous in \([-1,1]\).
		
		\item The function $\kernel$ has series representation of the form
		\[
		f_{\sphere[]}(r,s, \rho)
		= \sum_{n=0}^\infty \alpha^{\dims}_n(r,s) \uPnorm_n(\rho), \qquad \rho\in [-1,1],
		\]
		where \(\alpha^{\dims}_n: \domain^2 \to \complex\) are \pdkernels{} for all \(n\in \nat_0\)
		that are summable, i.e.  \(\sum_{n=0}^\infty \alpha^{\dims}_n(r,r) <
		\infty\) for all \(r\in \domain\).
	\end{enumerate}
	If $\kernel$ satisfies one of these equivalent conditions and \(\dims<\infty\)
	the \shortkernels{} \(\alpha^{\dims}_n\) moreover have the representation
	\[
	\alpha^{\dims}_n(r,s)
	= \frac{\scp{f_{\sphere[]}(r,s, \cdot)}{\uPnorm_n}_{\weight}}{\norm{\uPnorm_n}_{\weight}^2},
	\]
	with \(\weight(\rho) = (1-\rho^2)^{\frac{\dims-3}2}\) and norm \(\norm{\cdot}_{\weight}\) induced by
	\[
		\scp{f}{g}_{\weight} = \int_{-1}^{1} f(\rho) \conj{g(\rho)} \weight(\rho) d\rho.
\]
\end{theorem}

\begin{lemma}[Characterization of \(O(\dims)\),  {\citep[e.g.][Prop.\ 4.A.2]{benningDistributionalViewHigh2025}}]
	\label{lem: characterization of O(d)}
	For \(x_1, \dots, x_n \in \real^\dims\), \(y_1,\dots, y_n \in \real^\dims\) with \(\dims \in \nat\cup \{\infty\}\) the following 
	are equivalent
	\begin{enumerate}[label={(\alph*)}]
		\item\label{it: inner products equal} \(\langle x_i, x_j\rangle = \langle y_i, y_j\rangle\) for all \(i,j=1,\hdots, n\),
		\item\label{it: linear isometry exists} there exists \(\phi\in O(\dims)\) such that \(\phi(x_i) = y_i\) for all \(i=1,\hdots,n\).
	\end{enumerate}
\end{lemma}

The following observation is used for the extension of an isotropic continuous  \pdkernel{} on $(0,\infty) \times \sphere \simeq \real^\dims\setminus \{0\}$ to a \pdkernel{} on  $\real^\dims$.
\begin{lemma}[Continuous extension]
	\label{lem: continuous extension}
	Let \(\kernel\colon \domain \times \domain \to \complex\) be a continuous function
	on some topological space \(\domain\) and let $A\subset \domain$ be a
	dense subset of \(\domain\).
	If $\kernel$ is positive definite on $A$, then \(\kernel\) is positive definite on \(\domain\).
\end{lemma}
\begin{proof}
	Let \(x_0,\dots, x_m\in \domain\) and \(c_0,\dots, c_m \in \complex\)
	and choose for $i=0,\dots, m$ a sequence $(x_i^{(n)})_{n\in\nat} \subset  A$ such that $x_i^{(n)}\overset{n\to \infty}{\longrightarrow} x_i$.
	By continuity and positive definiteness of $\kernel$ on $A$ we have 
	\[
	\sum_{i,j=0}^m c_i \conj{c_j} \kernel(x_i, x_j)= \lim_{n\to \infty} \sum_{i,j=0}^m  c_i\conj{c_j} \kernel(x_i^{(n)}, x_j^{(n)})\ge 0.
	\qedhere
	\]
\end{proof}
\subsection{Proofs for the main part}\label{sec:all-proofs}

\subsubsection{Section \ref{sec: introduction}: Dimensional hierarchy}

\begin{proof}[Proof of Lemma \ref{lem: valid in all dimensions}]
	`\(\supseteq\)' follows from \eqref{eq: pd subset relation}. For
	`\(\subseteq\)' we need to prove that the representation \(\kappa\) of a \pdkernel{} $\kernel$ on $\cap_{\dims \in \nat} \PDset_{O(d)}(\real^\dims)$ induces a \pdkernel{} on \(\ell^2\). For this consider vectors \(x_1,\dots,x_n \in \ell^2\)
	and apply Lemma \ref{lem: characterization of O(d)} to obtain $U\in
	O(\infty)$ to map \(x_1,\dots,x_n\) into \(\real^n\subseteq \ell^2\). By positive
	definiteness in \(\PDset_{O(n)}(\real^n)\) the definition of positive definiteness
	\eqref{eq:positive_definite_definition} follows.
\end{proof}

\subsubsection{Section \ref{sec:positive_definite}: Characterization of isotropic kernels}

\begin{proof}[Proof of Theorem \ref{thm: isotropy characterization}]
\textbf{\ref{it: positive, definite isotropic} \(\Leftrightarrow\) \ref{it: pos definite, weak representation}:} 
Apart from the continuity of \(\kappa\), this equivalence has already been shown
in \citet[Appendix F]{benningRandomFunctionDescent2024}.
\\
``\(\Leftarrow\)'': This follows
directly from the fact that \(U\) is a linear isometry that preserves norms
and inner products. \\
``\(\Rightarrow\)'':
Let $M_\ge := \{(r,s,\gamma) \in M : r \ge s\}$.  Using two orthonormal
vectors \(e_1, e_2\), which exist due to \(\dims\ge 2\), we define
\emph{continuous} maps
$\varphi,\psi: M_\ge \to \real^\dims$ by
\[
\varphi(r,s,\gamma) := r e_1,\qquad
\psi(r,s,\gamma) := \begin{cases}
	\frac{\gamma}{r} e_1 + e_2 \sqrt{s^2 - \frac{\gamma^2}{r^2}} & r>0,\\
	0 & r=0.
\end{cases}
\]
By construction of $M_{\ge}$ all \((r,s,\gamma)\in M_\ge\) satisfy $r\ge s$ and $rs\ge \vert\gamma\vert$ and hence $\psi$ is continuous away from $(0,0,0)$, and continuity in $(0,0,0)$ follows directly from $\norm{\psi(r,s,\gamma)}=s$. 
The goal of the auxiliary functions $\varphi$ and $\psi$ is to create two vectors of
specified length and angle. Indeed, for \(\theta=(r,s,\gamma)\) we have
\begin{equation}
	\label{eq:preserved_directions}
	\|\varphi(\theta)\| = r, \quad \|\psi(\theta)\| = s, \quad \langle \varphi(\theta), \psi(\theta)\rangle = \gamma.
\end{equation}
The continuity requirement prevented the definition of \(\varphi\) and \(\psi\)
on the entirety of $M$. To deal with the cases $r<s$ we define the swap function
\(\tau(r,s,\gamma) := (s,r,\gamma)\).
With its help we finally define $\kappa\colon M\to \complex$ for $\theta = (r,s,\gamma)$ as
\begin{equation}
	\label{eq: definition of kappa}
	\kappa(r,s,\gamma) := 
	\begin{cases}
		\kernel\bigl(\varphi(\theta), \psi(\theta)\bigr) &  \text{ if } r\ge s \iff \theta \in M_\ge  \\
		\kernel\bigl(\psi(\tau(\theta)), \varphi(\tau(\theta)) \bigr)	& \text{ if } r\le s.
	\end{cases}
\end{equation}
Clearly, \(\kappa\) is continuous if it is well-defined in \(r=s\) since \(\varphi\), \(\psi\) and \(\tau\) are
continuous. For the case \(r=s\),
Lemma~\ref{lem: characterization of O(d)} and equation
\eqref{eq:preserved_directions} yield a linear isometry $U\in
O(\dims)$ such that $U\varphi (\theta) = \psi(\tau(\theta))$
and $U \psi(\theta) = \varphi(\tau(\theta))$ and thus
\[
\kernel\bigl(\varphi(\theta) , \psi(\theta)\bigr)
\overset{\text{isotropy}}= \kernel\bigl( U \varphi(\theta), U\psi(\theta) \bigr)
= \kernel\bigl(\psi(\tau(\theta)), \varphi(\tau(\theta))\bigr).
\]
It remains to show that $\kappa(\| x\|,\|y\|,\langle x,y\rangle)= \kernel
(x,y)$ holds for all $x,y\in \real^\dims$. Since the other case works analogously
we assume w.l.o.g.\@ $\|x\|\ge \|y\|$ and denote by $\theta :=(\|x\|,\|y\|,\langle x,y
\rangle)$ the input to \(\kappa\). Then by Lemma~\ref{lem: characterization of O(d)}
and \eqref{eq:preserved_directions} there again exists
$U\in  O(\dims)$  such that \(U x = \varphi(\theta)\) and \(Uy = \psi(\theta)\)
and we conclude
\begin{align*}
	\kernel(x,y)
	= \kernel(Ux, Uy)
	=  \kernel(\varphi(\theta), \psi(\theta))
	= \kappa(\theta).
\end{align*}

\textbf{\ref{it: positive, definite isotropic} \(\Rightarrow\) \ref{it: representation}:}
We begin by calculating the series representation of the kernel $\kernel$ on $\real^\dims\setminus\{0\}$.
Using the homeomorphism \(T\colon (0,\infty) \times \sphere\to
\real^d\setminus\{0\}\) with \(T(r,v) = rv\)
we define the continuous
\shortkernel{} $\kernel_{\sphere[]}\colon ((0,\infty) \times \sphere)^2 \to \complex $ by
\[
	\kernel_{\sphere[]} ((r,v),(s,w))
	\coloneq \kernel(rv,sw)
	\overset{\text{\ref{it: pos definite, weak representation}}}=
	\kappa\bigl(r, s, rs \langle v, w\rangle \bigr).
\]
With \(f_{\sphere[]}(r,s,\rho)\coloneq \kappa(r,s,rs\rho)\)
it is clear that Theorem~\ref{thm: guella XxS characterization} is
applicable and we obtain a representation in terms of \(f_{\sphere[]}\) for any
\(x,y\in \real^\dims\setminus\{0\}\), that is
\[
\kernel(x,y)
= f_{\sphere[]}\Bigl(\|x\|,\|y\|,\bigl\langle \tfrac{x}{\|x\|},\tfrac{y}{\|y\|}\bigr\rangle\Bigr)
= \sum_{n=0}^\infty \alpha_n^{\dims} (\|x\|,\|y\|) \uPnorm_n\Bigl(\bigl\langle \tfrac{x}{\|x\|},\tfrac{y}{\|y\|}\bigr\rangle\Bigr),
\]
where $\alpha_n^{\dims}$ are \pdkernels{} on
$(0,\infty)$ that are summable on the diagonal.
It remains to extend this characterization to \(0\) and
prove the claimed properties of the \(\alpha_n^\dims\).
To this end we investigate the cases $\dims<\infty$ and $\dims=\infty$
separately before proving the local uniform convergence of the partial sums of
\(\alpha_n^\dims\) as the last step of the proof.

\begin{enumerate}[wide=0pt]
	\item \textbf{Case \(\dims<\infty\)}:
	For all $n\in\nat$ and all $r,s\in (0,\infty)$ we have by Theorem~\ref{thm: guella XxS characterization} 
	\begin{equation}
		\label{eq: alpha formula}
		\alpha^{\dims}_n(r,s)
		= \frac{\scp{f_{\sphere[]}(r,s, \cdot)}{\uPnorm_n}_{\weight}}{\norm{\uPnorm_n}_{\weight}^2}
		= \frac{1}{\norm{\uPnorm_n}_{\weight}^2}
		\int_{-1}^1 \kappa(r,s,rs\rho) \uPnorm_n(\rho)\weight(\rho) d\rho.
	\end{equation}
	Equation \eqref{eq: alpha formula} is the claimed representation of the \(\alpha_n^\dims\) in \eqref{eq:alpha_formula}
	restricted to \((0,\infty)^2\). We use \eqref{eq: alpha formula} as
	definition to extend \(\alpha^{\dims}_n\) to a function on
	$[0,\infty)^2$. In particular, we get \eqref{eq:alpha_formula} by definition. To see
	that the extended \(\alpha_n^\dims\) are continuous, we use that
	\(\kappa\) is continuous by \ref{it: pos definite, weak
		representation} and apply dominated convergence to \eqref{eq: alpha formula} using the integrable weight function \(\weight(\rho) = (1-\rho^2)^{\frac{\dims-3}2}\) and the bound
	\[
	|\uPnorm_n(\rho)\weight(\rho)|
	\le \weight(\rho)
	\]
	which follows from \(|\uP_n(\rho)| \le \uP_n(1)\) \citep[e.g.][Eq. (2.16)]{reimerGegenbauerPolynomials2003}.
	This continuous extension of \(\alpha_n^\dims\) remains positive definite by Lemma~\ref{lem:
		continuous extension}. To see that this extension results in a representation
	of $\kernel$, observe that, by the orthogonality of
	the Gegenbauer polynomials with respect to the weight function \(\weight\)
	\citep[Theorem~2.3]{reimerGegenbauerPolynomials2003} and $\uPnorm_0 \equiv 1$
	we have for $r\ge 0$
	\[
	\alpha^{\dims}_n(r,0)
	= \frac{\bigl\langle \kappa(r,0, 0), \uPnorm_n\bigr\rangle_\weight}{\|\uPnorm_n\|_w^2}
	= \kappa(r,0, 0)\frac{\bigl\langle \uPnorm_0, \uPnorm_n\bigr\rangle_\weight}{\|\uPnorm_n\|_w^2}
= \kappa (r,0,0) \delta_{n0},
	\]
	where $\delta_{n0}$ denotes the Kronecker delta.
	This implies \ref{it: zero at the origin},
	i.e.\ the kernels \(\alpha_n^\dims\) are zero at the origin for all
	\(n\ge 1\), and we obtain the representation of the kernel in the origin
	in Equation \eqref{eq: kernel representation}, that is
	\[
		\kernel(x,0)
		= \kappa(\|x\|,0, 0 )
		= \alpha_0^\dims(\|x\|, 0).
	\]
	By an analogous argument the representation also holds for
	\(\kernel(0,x)\). The representation is also unique by
	orthogonality of the Gegenbauer polynomials, because if
	\(\tilde\alpha^{\dims}_n\) is another representation, then
	\begin{align*}
		\alpha_m^{\dims}(r,s)
		= \frac{\langle \kappa(r,s, rs\,\cdot), \uPnorm_m\rangle_\weight}{\|\uPnorm_m\|_\weight^2}
		&= \frac{\Bigl\langle \sum_{n=0}^\infty \tilde\alpha_n^{\dims}(r,s)\uPnorm_n(\cdot), \uPnorm_m(\cdot)\Bigr\rangle_\weight}{\|\uPnorm_m\|_\weight^2}
		\\
		&= \tilde\alpha_m^{\dims}(r,s).
	\end{align*}
	From the integral representation \eqref{eq:alpha_formula} we furthermore see that all $\alpha_n^{\dims}$  are real if $\kernel$ is real, as the Gegenbauer polynomials are real valued.
	The other direction follows directly from the series representation of $\kernel$. That is, $\kernel$ is real if and only if all $\alpha_n^{\dims}$ are real.
	
	\item \textbf{Case \(\dims=\infty\):}
	First note that the \(\alpha^{\infty}_n\) are \emph{continuous} \pdkernels{} on
	\((0,\infty)\)
	by \cite[Example~2.4]{guellaSchoenbergsTheoremPositive2019}.  For the
	extension to the origin we select two orthonormal vectors \(e_1,e_2\) and observe
	\begin{align*}
		\kernel(r e_1, r e_1) - \kernel(r e_1, r e_2)
		&= \sum_{n=0}^\infty \alpha^{\infty}_n(r,r)\langle e_1, e_1\rangle^n
		- \sum_{n=0}^\infty \alpha^{\infty}_n(r,r)\langle e_1, e_2\rangle^n
		\\
		&= \sum_{n=\red{1}}^\infty \alpha^{\infty}_n(r,r).
	\end{align*}
	By continuity of \(\kernel\) we obtain
	\[
	0 = \lim_{r\to 0}\kernel(r e_1, r e_1) - \kernel(r e_1, r e_2)
	= \lim_{r\to 0}\sum_{n=1}^\infty \alpha_n^{\infty}(r,r).
	\]
	Since \(\alpha^{\infty}_n(r,r)\ge 0\) holds for all $n\in\nat$ as
	the $\alpha^{\infty}_n$ are kernels, this implies
	\(\alpha^{\infty}_n(r,r) \to 0\) for all \(n\ge 1\).
	Assume \(r=0\) or \(s=0\) and 
	let \((r_l, s_l)\to(r,s)\) as \(l\to \infty\). Then the Cauchy-Schwarz inequality yields
	\[
	|\alpha^{\infty}_n(r_l, s_l) - \underbrace{\alpha^{\infty}_n(r,s)}_{=0}|
	\overset{\text{C.S.}}\le \sqrt{\alpha^{\infty}_n(r_l, r_l) \alpha^{\infty}_n(s_l, s_l)} \overset{l\to\infty} \longrightarrow 0.
	\]
	Hence, we can continuously extend \(\alpha^{\infty}_n\) to be zero in the origin \ref{it: zero at the origin}.
	Lemma~\ref{lem: continuous extension} implies for \(n\ge 1\) that this continuous extension of \pdkernels{} \(\alpha^{\infty}_n\)  is positive definite. 
	We continue with the extension in the case \(n=0\). For this we observe that
	\[
		\kernel(r e_1, s e_2)
		= \sum_{n=0}^\infty \alpha^{\infty}_n(r,s) \langle e_1, e_2\rangle^n
		= \alpha^{\infty}_0(r,s),
	\]
	using \(\scp{e_1}{e_2} = 0\) for \(r,s>0\) and for $r=0$ or $s=0$ that the kernels $\alpha_n^\infty(r,s)$ are zero for $n\ge 1$.
	This yields the continuous extension \(\alpha^{\infty}_0(r,s) :=
	\kernel(re_1,se_2)\), which remains positive definite by Lemma~\ref{lem: continuous extension}.
	What is left to prove is
	that our continuous extension of the \shortkernels{}
	\(\alpha^{\infty}_n\) does in fact lead to a representation 
	of \(\kernel\) and the uniqueness of the representation. 
	For the former, let $x,y\in \real^\dims$ and w.l.o.g.\@ assume $x=0$ to
	obtain by isotropy
	\[
	K(x,y)
	= K(0,\|y \| e_2 )
	= \alpha_0^{\infty}(0,\|y\| ).
\]
	For uniqueness of the \(\alpha_n^{\infty}\) observe that
	\(
	\kappa(r,s, rs\rho) = \sum_{n=0}^\infty \alpha_n^{\infty}(r,s) \rho^n
	\)
	is clearly an analytic function in \(\rho\) defined on \([-1,1]\) for
	all \(r,s\) and the coefficients \(\alpha_n^{\infty}(r,s)\) thereby unique.
	This also shows that all kernels $\alpha^{\infty}_n$ are real valued if and
	only if $\kernel$ is real valued. This is because the coefficients of
	a real-valued analytic function are the real valued derivatives of said function.
\end{enumerate}

For the local uniform convergence of the partial sums \(\alpha_n^\dims\), \ref{it: locally uniformly convergent},
we use the continuity of \(\kernel\) and the continuity of the \(\alpha_n^\dims\) we obtained 
in the previous steps. Indeed, since the \(\alpha_n^\dims\) are continuous,
the partial sums
\begin{equation}
	\label{eq: partial sums}
	S_N(r) \coloneq \sum_{n=0}^N \alpha_n^\dims(r,r)
\end{equation}
are continuous in \(r\) and converges pointwise to the continuous function
\[
	S(r) \coloneq \kernel(re_1,re_1) = \sum_{n=0}^\infty \alpha_n^\dims(r,r)\uPnorm_n(1)
\]
using \(\uPnorm_n(1) = 1\). Since \(\alpha_n^\dims(r,r)> 0\)
due to positive definiteness of the \(\alpha_n^\dims\), the \(S_N\) are monotonically increasing,
which means we may apply Dini's theorem \citep[e.g.][Thm.~7.13]{rudinPrinciplesMathematicalAnalysis2008} to conclude that the convergence is uniform on any
compact interval \([0,R]\) for any \(R>0\), that is
\[
	\sup_{r\in [0,R]} \abs[\Big]{\sum_{n=N}^\infty \alpha_n^{\dims}(r,r)}
	= \sup_{r\in [0,R]} \abs[\Big]{\sum_{n=0}^N \alpha_n^{\dims}(r,r) - \kernel(rv, rv)} \to 0
	\quad \text{for } N\to \infty.
\]
In particular for any \(r\in [0,\infty)\) we can find a neighborhood \(U_r\) of \(r\) contained
in \([0,R)\) for some \(R>0\) such that the convergence is uniform on \(U_r\).

\textbf{\ref{it: representation} \(\Rightarrow\) \ref{it: positive, definite isotropic}:} 
For the continuity of \(\kernel\) define the partial sums
\[
	\kernel_N(x,y) \coloneq \sum_{n=0}^N \alpha_n^{\dims}(\norm{x}, \norm{y}) \uPnorm_n\bigl(\scp[\big]{\tfrac{x}{\norm{x}}}{\tfrac{y}{\norm{y}}}\bigr).
\]
To see that the \(\kernel_N\) are continuous, consider the individual summands.
For \(n=0\) we have continuity by continuity of \(\alpha_0^\dims\) and
\(\uPnorm_0\equiv 1\). And for \(n\ge 1\) use the continuity of \(\alpha_n^\dims\) and \(\uPnorm_n\) on \([-1,1]\) and
the fact that the \(\alpha_n^\dims\) are zero at the origin and
\(\abs{\uPnorm_n(\rho)} \le 1\) for all \(\rho\in [-1,1]\).
Moreover, the functions \(\kernel_N\) converge locally uniformly to \(\kernel\) by the local uniform convergence of the partial
sums on the diagonal due to
\begin{align*}
	\abs{\kernel(x,y) - \kernel_N(x,y)}
	&\le \sum_{n=N}^\infty \abs[\big]{\alpha_n^{\dims}(\norm{x}, \norm{y})} \abs[\Big]{\uPnorm_n\bigl(\scp[\big]{\tfrac{x}{\norm{x}}}{\tfrac{y}{\norm{y}}}\bigr)}
	\\
	&\le \sum_{n=N}^\infty \sqrt{\alpha_n^{\dims}(\norm{x}, \norm{x})\alpha_n^{\dims}(\norm{y}, \norm{y})}
	\\
	&\le \frac12\sum_{n=N}^\infty \alpha_n^{\dims}(\norm{x}, \norm{x}) + \alpha_n^{\dims}(\norm{y}, \norm{y})
	\to 0 \qquad (N\to \infty),
\end{align*}
where the convergence to zero is uniform in \(x\) and \(y\) on some neighborhoods \(U_{\norm{x}}\) and \(U_{\norm{y}}\) of \(\norm{x}\) and \(\norm{y}\)
by \ref{it: locally uniformly convergent}. But a sequence of continuous, locally uniformly convergent
functions \(\kernel_N\) converges to a continuous limit \(\kernel\) \citep[e.g.][Thm.~7.12]{rudinPrinciplesMathematicalAnalysis2008}.
To apply the cited result the reader should restrict the domain to the local neighborhoods to
turn the local uniform convergence into a uniform convergence.

Isotropy is obvious and it remains to show
positive definiteness of $\kernel$. We split the series representation into two parts
\begin{equation}
	\label{eq: kernel decomposition}
	\kernel(x,y)
	= \underbrace{\alpha^{\dims}_0(\|x\|, \|y\|)}_{\defEqual\kernel_0(x,y)}
	+ \underbrace{\sum_{n=1}^\infty \alpha^{\dims}_n(\|x\|, \|y\|)\uPnorm_n\bigl(\langle \tfrac{x}{\|x\|}, \tfrac{y}{\|y\|}\rangle\bigr)}_{\defEqual\kernel_{>0}(x,y)},
\end{equation}
where we use \(\uPnorm_0 \equiv 1\).  
Let \(m\in \nat\), \(x_1,\dots,x_m\in \real^\dims\) and \(c_1,\dots,c_m\in \complex\),
and assume without loss of generality $x_i=0$ if and only if $i=1$. We can do so with the choice \(c_1=0\)
if none of the \(x_i\) are zero. Since
$\alpha_n^{\dims}(0,r)=0$ for $n\ge 1$ we obtain
\[
\sum_{i,j=1}^m c_i \conj{c_j} \kernel(x_i,x_j)
= \sum_{i,j=1}^m c_i\conj{c_j} \kernel_0(x_i,x_j) + \sum_{i,j=2}^m c_i\conj{c_j} \kernel_{>0}(x_i,x_j).
\]
Since $\alpha_0^{d}$ is a \shortkernel{} on $[0,\infty)$ we get that the first
summand is non-negative. The second summand is non-negative, since it only
contains $x_i\in \real^\dims \setminus\{0\}$ and the restriction
$\kernel_{>0}\vert_{(\real^\dims\setminus \{0\})^2}$ is positive definite
using the homeomorphism $T:(0,\infty) \times \sphere\to \real^\dims
\setminus \{0\},\; (r,v)\mapsto rv$ and Theorem~\ref{thm: guella XxS
	characterization} characterizing \shortkernels{} on \((0,\infty)\times \sphere\).

\end{proof}

\subsubsection{Section \ref{sec: strict positive definite}: Characterization of strict positive definite kernels}

\begin{proof}[Proof of Theorem \ref{thm:strictly_pos_def}]
	Recall the decomposition 
	\eqref{eq: kernel decomposition}, i.e.\
	\begin{equation}
		\label{eq: recall kernel decomposition}	
		\kernel(x,y)
		= \underbrace{\alpha^{\dims}_0(\|x\|, \|y\|)}_{=\kernel_0(x,y)}
		+ \underbrace{\sum_{n=1}^\infty \alpha^{\dims}_n(\|x\|, \|y\|) \uPnorm_n(\langle \tfrac{x}{\|x\|}, \tfrac{y}{\|y\|}\rangle)}_{=\kernel_{>0}(x,y)},
	\end{equation}
	with \pdkernels{} $\kernel_0$ and $\kernel_{>0}$.  

	``\ref{it: infinitely even and odd elements in set}, \ref{it: even and odd kernels strict pos. def.}, \ref{it: infinitely many even and odd n strict pos. def.} \(\Rightarrow\) \ref{it: strict positive definite}'':
	For strict positive definiteness of \(\kernel\) we need to show
	\[
		\sum_{i,j=1}^m c_i \conj{c_j} \kernel(x_i,x_j)
	\]
	is non-zero for any \(m\in \nat\), \(c_1, \dots, c_m\in
	\complex\setminus\{0\}\) and distinct  \(x_1,\dots, x_m\in \real^\dims\).
	If none of the \(x_i\) are zero the result follows directly from 
	Theorem 3.5, 3.7 and 3.9 of \citet{guellaSchoenbergsTheoremPositive2019} with the usual
	homeomorphism to \((0,\infty)\times \sphere\) as the conditions \ref{it:
	infinitely even and odd elements in set}, \ref{it: even and odd kernels
	strict pos. def.} and \ref{it: infinitely many even and odd n strict pos.
	def.} imply the conditions of Theorem 3.5, 3.7 and 3.9 respectively.
	If at least one \(x_i\) is zero, we can assume \(x_1=0\) w.l.o.g.\@, then we
	have
	\begin{equation}
		\label{eq: decomposition}
		\sum_{i,j=1}^m c_i \conj{c_j} \kernel(x_i,x_j)
		= \sum_{i,j=1}^m c_i\conj{c_j} \kernel_0(x_i,x_j)
		+ \sum_{i,j=2}^m c_i\conj{c_j} \kernel_{> 0}(x_i,x_j)
\end{equation}
	using that \(\kernel_{>0}(x,0) = \kernel_{>0}(0,x) = 0\) for all \(x\in \real^\dims\)
	due to \eqref{eq: zero in origin} in Theorem~\ref{thm: isotropy
	characterization}. Since both summands are non-negative by positive-definiteness it is sufficient if
	either one is non-zero.
	
	In the case \(m>1\) we observe that the conditions \ref{it: infinitely even and odd elements in set},
	\ref{it: even and odd kernels strict pos. def.} and \ref{it: infinitely many
	even and odd n strict pos.  def.} imply that \(\kernel_{>0}\) mapped to
	\((0,\infty)\times \sphere\) satisfies the conditions of Theorem 3.5, 3.7
	and 3.9 of \citet{guellaSchoenbergsTheoremPositive2019} and hence is
	strictly positive definite on \(\real^\dims\setminus\{0\}\). This implies
	that the second summand in \eqref{eq: decomposition} is strictly greater
	than zero.

	And in the case $m=1$ the first term is strictly greater than zero, since
	\[
	\sum_{i,j=1}^m c_i\conj{c_j} \kernel_0(x_i,x_j)
	= |c_1|^2 \alpha_0^{\dims}(0,0) >0.
	\]
	
	``\ref{it: strict positive definite} \(\Rightarrow\) \ref{it: infinitely even and odd elements in set}, \ref{it: even and odd kernels strict pos. def.}'':
	If \(\kernel\) is strictly positive definite on \(\real^\dims\), then
	we have \(\alpha^{\dims}_0(0,0)>0\) by the selection of \(m=1\), \(c=1\),
	\(x_1=0\) cf.~\eqref{eq: recall kernel decomposition}. 
	Mapping strict positive definiteness of \(\kernel\) on \(\real^\dims\setminus\{0\}\)
	into strict positive definiteness on \((0,\infty)\times \sphere\), the
	rest of the assertions in \ref{it: infinitely even and odd elements in set} and \ref{it:
		even and odd kernels strict pos. def.} follow from Theorem 3.5 and 3.7 of
	\citet{guellaSchoenbergsTheoremPositive2019}.
	
	The final remark about strict positive definiteness on \(\real^\dims\setminus\{0\}\)
	follows similarly from the homeomorphism to \((0,\infty)\times \sphere\) and
	the results of \citet{guellaSchoenbergsTheoremPositive2019}.
\end{proof}

\begin{proof}[Proof of Theorem \ref{thm:strictly_pos_def_dim2}]
	The proof of this Theorem is identical to the proof of Theorem \ref{thm:strictly_pos_def}	
	with the uses of Theorems 3.5, 3.7 and 3.9 by
	\citet{guellaSchoenbergsTheoremPositive2019} replaced by Theorems 3.6, 3.8
	and 3.9 respectively.
\end{proof}

\begin{remark}[Universal kernels]
	Universal kernels are kernels whose associated reproducing kernel Hilbert
	space is dense in the space of compactly supported continuous functions with respect to the
	supremum norm \citep[Def.~2 and Thm.~1 in Sec.~4]{carmeliVectorValuedReproducing2010}. This property provides universal learning guarantees for kernel
	methods using such kernels.
	For future work characterizing \emph{universal} isotropic \shortkernels{} we point out
	their characterization on the sphere by \citet[Theorem
	10]{micchelliUniversalKernels2006}.  It might be possible to extend this
	result to $(0,\infty) \times \sphere$ analogous to the work of
	\citet{guellaSchoenbergsTheoremPositive2019} for (strictly) \pdkernels{},
	and if that succeeds one can use the isomorphy of $(0,\infty)\times \sphere$
	and $\real^\dims \setminus \{0\}$ in a fashion similar to our proof of
	Theorem \ref{thm: isotropy characterization} and Theorem
	\ref{thm:strictly_pos_def} to obtain the desired characterization.
\end{remark}

\subsection{Alternative proof of the characterization of continuous dot product kernels \texorpdfstring{on \(\ell^2\)}{}}
The following lemma, which can be found in \citep{pinkusStrictlyPositiveDefinite2004} and the references
therein, classifies which functions $\kappa\colon\real\to\real$ define (strictly)
positive definite dot product kernels valid in all dimensions, or
equivalently in $\ell^2$.
Using Theorem \ref{thm: isotropy characterization} we provide an
alternative proof. 

\begin{lemma}[Classification of continuous dot product kernels on $\ell^2$]
	\label{lem: classification of dot product kernels infinite dim}
	Let $\kernel(x,y) := \kappa(\scp xy)$ be a continuous positive definite scalar product kernel on $\ell^2$. Then 
	\begin{align}
		\label{eq:scalar_kernel_classification}
		\kernel(x,y) = \sum_{n=0}^\infty a_n \scp xy^n
	\end{align}
	holds, where $0\le a_n$ for all $n$ and  $(a_n t^n )_{n\ge 0} \in \ell^1$ holds for all $t\in\real$. On the other hand, each such sequence $(a_n)_{n\ge 0}$ yields a continuous scalar product kernel on $\ell^2$.
	
	Furthermore a kernel  $\kernel(x,y)= \sum_{n=0}^\infty a_n \scp xy^n$ is strictly positive definite on $\ell^2$ if and only if both $a_n>0$ holds for infinitely many even and infinitely many odd indices $n$, as well as $a_0>0$.
\end{lemma}
\begin{proof}
	If $\kernel(x,y) := \kappa(\scp xy)$ defines a positive definite kernel on $\ell^2$ we have 
	by Theorem \ref{thm: isotropy characterization}
	\[
		\kernel(x,y)
		= \sum_{n=0}^{\infty}\alpha_n(\norm x,\norm y)
		\bigl(\bigl\langle \tfrac{x}{\|x\|}, \tfrac{y}{\|y\|}\bigr\rangle\bigr)^n
		=\sum_{n=0}^{\infty}
		\underbrace{\frac{\alpha_n(\norm x,\norm y)}{(\norm x\norm y)^n}}_{
			=: \tilde\alpha_n(\|x\|, \|y\|)
		}
		\scp xy^n.
	\]
	Since the \(\alpha_n\) are continuous we may restrict ourselves to \(x,y\in \ell^2\setminus\{0\}\)
	where \(\tilde \alpha_n\) is well defined. To prove \(\tilde \alpha_n\) is constant
	for every \(n\) we need to show for any \(r_1,s_1 > 0\) and \(r_2, s_2>0\)
	that \(\tilde \alpha_n(r_1,s_1) = \tilde \alpha_n(r_2,s_2)\). To this end
	we define for \(i\in \{1,2\}\)
	\[
		x_i(t) = r_i e_1
		\qquad
		y_i(t) = \frac{t}{r_i}e_1 + \sqrt{s_i^2 - (t/r_i)^2} e_2
	\]
	for all \(t\in [-\epsilon, \epsilon]\) with \(\epsilon = \min\{r_1s_1, r_2s_2\}>0\).
	Observe that these satisfy
	\[
		\scp{x_i(t)}{y_i(t)}=t,
		\qquad
		\norm{x_i(t)} = r_i
		\quad\text{and}\quad
		\norm{y_i(t)}=s_i.
	\]
	Consequently we have for all \(t\in [-\epsilon, \epsilon]\)
	\begin{align*}
		0 &= \kappa(t) - \kappa(t)
		\\
		&= \kernel(x_1(t), y_1(t)) - \kernel(x_2(t), y_2(t))
		\\
		&= \sum_{n=0}^\infty \bigl(\tilde \alpha_n(r_1, s_1) - \tilde \alpha_n(r_2,s_2)\bigr) t^n.
	\end{align*}
	This implies \(\tilde \alpha_n(r_1, s_1) = \tilde \alpha_n(r_2,s_2)\) because the
	coefficients of a power series are unique. And since the \(r_i\) and \(s_i\) were arbitrary
	the \(\tilde \alpha_n\) are constant, say \(\tilde \alpha_n \equiv a_n
	\in \real\) on \((0, \infty)^2\).
	By Theorem \ref{thm: isotropy characterization}, the $\alpha_n$ are continuous
	kernels on $[0,\infty)^2$ that are zero at the origin for \(n\ge 1\) and we thus obtain
	\[
		\alpha_n(\norm x,\norm y) = a_n \cdot (\norm x\norm y)^n,
	\]
	where the \(a_n\) must be non-negative for \(\alpha_n\) to be positive
	definite.
	But this implies \eqref{eq:scalar_kernel_classification} and we observe that
	$(a_nt^n)_{n\in \nat_0} \in \ell^1$ holds for all $t\ge 0$ as
	\eqref{eq:scalar_kernel_classification} converges. 
	
	For the converse statement observe that
	each non-negative series $a_n$ with $(t^na_n)_{n\in \nat_0} \in \ell^1$ for all $t\in
	\real$  defines a scalar product kernel on $\ell^2$ by invoking Theorem
	\ref{thm: isotropy characterization}. With this characterization of scalar
	product kernels, the claims for strict positive definiteness follow directly
	from Theorem \ref{thm:strictly_pos_def}.
\end{proof}

\subsection{Technical lemma for Gaussian convergence}

The following lemma clarifies technical convergence questions in the proof of
Theorem \ref{thm: representation isotropic random functions}. It is essentially
a direct consequence of classical martingale theory \citep[e.g.][Theorem
11.10]{klenkeProbabilityTheoryComprehensive2014}.

\begin{lemma}
    \label{lem: series convergence}
    Let \((X_n)_{n\in \nat}\) be a sequence of independent, centered Gaussian 
    random vectors in \(\real^\dims\) with \(\sum_{n\in \nat} \E[\norm{X_n}^2] < \infty\).
    Then the sequence of partial sums \(S_N \coloneq \sum_{n=1}^N X_n\) converges almost surely
    and in \(L^2(\Omega)\) against a centered Gaussian vector \(S_\infty\). Moreover, the covariance
    matrix of \(S_\infty\) is given by the sum of covariance matrices, that is
    \[
        \Cov(S_\infty) = \sum_{n\in \nat} \Cov(X_n).
    \]
\end{lemma}
\begin{proof}
	Since the Gaussian random vectors $(X_n)_{n\in\nat}$ are centered and
	independent, each component $S_N^i$ of the partial sums \(S_N\) is a
	martingale in \(N\) with respect to the natural filtration satisfying
	\begin{align*}
		\E[(S_n^i)^2] = \sum_{k=1}^n \E[(X_k^i)^2] \le  \sum_{k=1}^\infty \E[\norm{X_k}^2]< \infty.
	\end{align*}
	Consequently, the sequence is an \(L^2\)-bounded martingale, by
	\citep[Theorem 11.10]{klenkeProbabilityTheoryComprehensive2014}  there exists
	a limit \(S_\infty^i\in L^2(\Omega)\) such that $S_N^i\to S_\infty^i$ in $L_2(\Omega)$ and almost
	surely. This also implies that the vector \(S_N\) converges in \(L^2(\Omega)\)
	and almost surely to the vector \(S_\infty\) with components \(S_\infty^i\).
	By Lévy's continuity theorem $S_\infty$ is a centered Gaussian random vector, and $L_2(\Omega)$ convergence yields
	\[
        \Cov(S_\infty)_{ij}
        = \lim_{N\to \infty} \Cov(S_N^{i}, S_N^{j})
        = \lim_{N\to \infty}\sum_{n=1}^N \Cov(X_n)_{ij}.
		\qedhere
    \]
\end{proof}
 	\section{Recapitulation of spherical harmonics}
\label{sec: spherical harmonics}

In this section we recapitulate some facts about spherical harmonic analysis.
In particular, we elaborate on Gegenbauer polynomials, which appear in our main
results (e.g.\ Theorem \ref{thm: isotropy characterization}). This may help the
reader gain a better intuition why Gegenbauer polynomials appear in our
results. The second reason for this recapitulation is that it forms a natural
basis for the recapitulation of Müller's definition
of a Bessel function \citep{mullerAnalysisSphericalSymmetries1998}, which differs from the standard definition (Section
\ref{sec: bessel functions}). This is necessary, since we are using Müller's results to get a representation
of the characteristic function of the uniform distribution
(Lemma~\ref{lem: characteristic function of uniform distribution on sphere}). And this
representation is key in the characterization of the finite dimensional stationary isotropic
kernels in Theorem~\ref{thm: stationary isotropic kernels}. We therefore
provide a translation of Müller's definition of Bessel functions to the
standard definition. Finally, we discuss the dimension of the space of spherical
harmonics, which is needed for the complexity analysis of the Gaussian random
function simulation in Section~\ref{sec: GRF representation}.
\begin{definition}[Homogeneous polynomial spaces]
    Let \(\homP_n\) be the space of  \textbf{\(\dims\)-variate, \(n\)-degree homogeneous polynomials}
    mapping \(\real^\dims\) to \(\complex\)
    \[
        \homP_n \coloneq \set[\big]{p \in \complex[x_1, \ldots, x_\dims]: p(\lambda x) = \lambda^n p(x),\; \forall \lambda \in \real, x\in \real^\dims}.
    \]
    Let \(\homHarmonic_n\) be the corresponding subspace of \textbf{harmonic} homogeneous polynomials
    \[
        \homHarmonic_n \coloneq \set[\big]{p \in \homP_n: \Delta p = 0},
        \qquad\text{where} \quad \Delta \coloneq \sum_{i=1}^\dims \frac{\partial^2}{\partial x_i^2}.
    \]
    Due to their defining property, homogeneous polynomials are uniquely
    determined by their restriction to the unit sphere \(\homPSphere_n \coloneq \homP_n \restriction_{\sphere}\).
    This space of functions \(\homPSphere_n\) is a subspace of \(L^2(\sphere)\) with the inner product
    \begin{equation}
        \label{eq: inner product on sphere}    
        \scp{f}{g}_{\sphere} \coloneq \frac1{\vol{\sphere}}\int_{\sphere} f(x) \conj{g(x)} \sigma(dx),
    \end{equation}
    where \(\sigma\) is the surface measure and \(\vol{\sphere}\) the area of
    the sphere.  The \textbf{spherical harmonics}, are the harmonic homogeneous
    polynomials restricted to the sphere \(\sphHarmonic_n \coloneq
    \homHarmonic_n \restriction_{\sphere}\).
\end{definition}

\citet{mullerAnalysisSphericalSymmetries1998} explains that the homogeneous,
harmonic polynomials \(\homHarmonic_n\) are natural spaces for the analysis of rotation, since:
\begin{enumerate}
    \item They are \textbf{invariant} with respect to rotations from \(O(\dims)\).
    That is, for any \(p\in \homHarmonic_n\) and \(U\in O(\dims)\),
    \(p_U(x) \coloneq p(Ux)\) is still a member of \(\homHarmonic_n\)
    \citep[§2.17]{mullerAnalysisSphericalSymmetries1998}.
    
    \item The spaces \(\homHarmonic_n\) are \textbf{irreducible}, that is, there are no
    invariant non-trivial subspaces \citep[§2, Thm.~3]{mullerAnalysisSphericalSymmetries1998}.
\end{enumerate}
The homogeneous polynomials \(\homP_n\) are also invariant, but due to the
invariant subspace \(\homHarmonic_n\) clearly not irreducible. These ``primitive''
spaces \(\homHarmonic_n\) decompose \(L^2(\sphere)\) in the following way:

\begin{theorem}[Decompositions]
    The following statements hold
    \begin{enumerate}[label=(\roman*),noitemsep]
        \item\label{it: orthogonal decomposition of L^2}
        The spherical harmonics are an orthogonal decomposition of
        \(L^2(\sphere)\), that is
        \[
            L^2(\sphere) = \bigoplus_{n=0}^\infty \sphHarmonic_n
            \qquad\text{with}\qquad
            \sphHarmonic_n \perp \sphHarmonic_m \quad \text{for } n\neq m.
        \]
        \item\label{it: decomposition of homogeneous polynomials}
        The homogeneous polynomial spaces decompose into
        \[
            \homP_n = \bigoplus_{k=0}^{\floor{n/2}} \norm{x}^{2k} \homHarmonic_{n-2k},
        \]
        that is for all \(p\in \homP_n\) there exist unique \(h_{n-2k} \in \homHarmonic_{n-2k}\) such that
        \[
            p(x) = \sum_{k=0}^{\floor{n/2}} \norm{x}^{2k} h_{n-2k}(x) \qquad \forall x\in \real^\dims.
        \]
        On the sphere \(\norm{x} = 1\) implies \(\homPSphere_n = \bigoplus_{k=0}^{\floor{n/2}} \sphHarmonic_{n-2k}\).
    \end{enumerate} 
\end{theorem}
\begin{proof}
    The orthogonality in \ref{it: orthogonal decomposition of L^2} is §2, Thm.~5 in \citep{mullerAnalysisSphericalSymmetries1998}
    or Theorem~1.1.2 in \citep{daiApproximationTheoryHarmonic2013}.
    \ref{it: decomposition of homogeneous polynomials} is Theorem~1.1.3 in \citep{daiApproximationTheoryHarmonic2013}.
    The decomposition part in \ref{it: orthogonal decomposition of L^2}
    now follows from \ref{it: decomposition of homogeneous polynomials} and the fact
    that all polynomials are sums of homogeneous polynomials, because the
    polynomials are dense in the continuous functions with respect to
    \(L^\infty(\sphere)\) by Stone-Weierstrass,
    see e.g.\ \citep[][Thm.~4.51]{follandRealAnalysisModern1999}, and thereby
    dense with respect to \(L^2\). Finally, the continuous functions are dense
    in \(L^2(\sphere)\) by \citep[][Prop.~7.9]{follandRealAnalysisModern1999}.
\end{proof}

Let \(N_n^\dims \coloneq \dim \sphHarmonic_n\) be the dimension of the space of spherical harmonics
with an orthonormal basis \(\set{Y_{n,1}, \ldots, Y_{n,N_n^\dims}}\)
of \(\sphHarmonic_n\), produced using Gram-Schmidt or explicit formulas, see
e.g.\ \citep[][Thm.~1.5.1]{daiApproximationTheoryHarmonic2013}. 
Then clearly
\begin{equation}
    \label{eq: zonal kernel}    
    \zonalK_n(x,y) \coloneq \sum_{k=1}^{N_n^\dims} Y_{n,k}(x) \conj{Y_{n,k}(y)}
\end{equation}
is a reproducing kernel for \(\sphHarmonic_n\) and thus \(\sphHarmonic_n\) is a
reproducing kernel Hilbert space. 
Since the reproducing kernel of an RKHS is unique,
\(\zonalK_n\) does not depend on the choice of the orthonormal basis. In particular
\(\zonalK_n\) is real valued since we may choose a real valued
orthonormal basis.
Moreover
\((Y_{n,1}(Q\cdot), \ldots, Y_{n,N_n^\dims}(Q\cdot))\) is another orthonormal basis for \(\sphHarmonic_n\),
because the space is invariant under rotations \citep[Theorem~1.1.7]{daiApproximationTheoryHarmonic2013}.
Therefore we have for all \(x,y\in \sphere, Q\in O(\dims)\) the identity \(\zonalK_n(Qx, Qy) = \zonalK_n(x,y)\) which yields
\[
    \zonalK_n(x,y) = \zonalK_n(\scp{x}{y}).
\]
Consequently, we may treat $\zonalK_n(t) , t\in [-1,1]$ like a univariate function
and \(\zonalK_n(x,x) = \zonalK_n(1)\) is constant for \(x\in\sphere\) with value
\begin{equation}
    \label{eq: scaling zonal kernel}    
    \zonalK_n(1) = \scp{\zonalK_n(x,x)}{1}_{\sphere}
    = \sum_{k=1}^{N_n^\dims} \scp{Y_{n,k}}{Y_{n,k}}_{\sphere}
    = N_n^\dims.
\end{equation}
Observe that by its definition in \eqref{eq: zonal kernel}, \(\zonalK_n(t)\) is a
polynomial in \(t\). 

Recall, that Gegenbauer polynomials are classically defined as the orthogonal
polynomials with respect to the weight \(\weight(t) = (1-t^2)^{\lambda -\frac{1}2}\)
\citep[e.g.][]{reimerGegenbauerPolynomials2003}.  The following ``addition
theorem''
\citep[e.g.][]{mullerAnalysisSphericalSymmetries1998,daiApproximationTheoryHarmonic2013}
for Gegenbauer polynomials may be viewed as an alternative definition
that is valid for \(\lambda = \frac{\dims-2}2\) with integer \(\dims \ge 2\).
This result may provide some intuition why Gegenbauer polynomials appear in our
characterization.

\begin{theorem}[Addition Theorem]
    \label{thm: addition theorem}
    The Gegenbauer polynomials are the reproducing kernels of \(\sphHarmonic_n\)
    up to a constant factor. Specifically for \(\lambda = \frac{\dims-2}{2}\)
    \[
        \uPnorm_n(\scp{x}{y})
        = \frac{\zonalK_n(x,y)}{N_n^\dims}
        = \frac1{N_n^\dims}\sum_{k=1}^{N_n^\dims} Y_{n,k}(x) \conj{Y_{n,k}(y)}.
    \]
\end{theorem}
The proof of this important theorem \citep[e.g.][]{mullerAnalysisSphericalSymmetries1998,daiApproximationTheoryHarmonic2013} is a simple verification of the defining orthogonality relations of the Gegenbauer polynomials
using the Funk-Hecke formula.
\begin{lemma}[Funk-Hecke {\citep[Lemma~A.5.2]{daiApproximationTheoryHarmonic2013}}] For \(x\in \real^\dims\) and \(f\colon \real\to \real\)
    \label{lem: Funk-Hecke}
    \[
        \int_{\sphere} f(\scp{x}{y}) \sigma(dy) = \vol{\sphere[\dims-2]}\int_{-1}^1 f(\norm{x}t) (1-t^2)^{\frac{\dims-3}{2}} dt.
    \]
\end{lemma}

Since the Addition theorem is where the Gegenbauer polynomials make their entrance, we
provide its short proof below for the reader's convenience.

\begin{proof}[Proof of Theorem~\ref{thm: addition theorem}]
Since \(\zonalK_n\) is real valued, the Funk-Hecke formula together with the
reproducing property of \(\zonalK_n\) implies the defining orthogonality relation
of Gegenbauer polynomials:
\[\begin{aligned}
    \int_{-1}^1 \zonalK_n(t) \zonalK_m(t) \underbrace{(1-t^2)^{\frac{\dims-3}{2}}}_{\weight(t)} dt
    \overset{\text{Lem.~\ref{lem: Funk-Hecke}}}&= \frac1{\vol{\sphere[\dims-2]}}\int_{\sphere} \zonalK_n(\scp{x}{y}) \zonalK_m(\scp{x}{y}) \sigma(dy)
    \\[-2ex]
    &= \frac{\vol{\sphere}}{\vol{\sphere[\dims-2]}} \scp[\big]{\zonalK_n(x,\cdot)}{\zonalK_m(x,\cdot)}_{\sphere}
    \\
    \overset{\sphHarmonic_n \perp \sphHarmonic_m}&=
    \frac{\vol{\sphere}}{\vol{\sphere[\dims-2]}}
    \underbrace{\zonalK_n(x,x)}_{=N_n^\dims} \delta_{n,m}.
\end{aligned}
\]
Since \(\zonalK_n\) is a polynomial of degree \(n\) by definition \eqref{eq: zonal kernel},
since \(\sphHarmonic_n\) is of degree \(n\) and since it satisfies the orthogonality
relation of the Gegenbauer polynomials, it must be a constant multiple by uniqueness.
The normalization \(\uPnorm_n(1)=1\) follows from \eqref{eq: scaling zonal kernel}.
\end{proof}

\subsection{Bessel functions}
\label{sec: bessel functions}

For the characterization of the stationary isotropic kernels in finite dimensions
in Theorem~\ref{thm: stationary isotropic kernels} we referenced a Bessel
functions representation of the characteristic function of the uniform
distribution on the sphere (Lemma~\ref{lem: characteristic function of uniform
distribution on sphere}). This Lemma~\ref{lem: characteristic function of uniform
distribution on sphere} is essentially Lemma 2 in §22 of
\citep{mullerAnalysisSphericalSymmetries1998}. However, \citet{mullerAnalysisSphericalSymmetries1998}
uses a custom definition of Bessel functions. We therefore explain in this section how his
definition translates to the classical definition of Bessel functions.

Up to constants \citet[§22, Def.~1]{mullerAnalysisSphericalSymmetries1998} \emph{defines}
a Bessel function to be the projection of \(e^{ir\cdot}\) onto the span
of the Gegenbauer polynomial \(\uPnorm_n\), specifically
\[
    \besselProj_n(\dims; r)
    \coloneq i^{-n} \frac{\vol{\sphere[\dims-2]}}{\vol{\sphere}}\int_{-1}^1 e^{irt} \uPnorm_n(t) \weight(t) dt
    \qquad \lambda = \tfrac{\dims-2}{2}.
\]
We note that his notation is \(P_n(t) = P_n(\dims; t)\) for \(\uPnorm_n(t)\)
as can be seen in \citep[(§8.8), (§9.1) and (§9.3)]{mullerAnalysisSphericalSymmetries1998}.
However, this definition of the Bessel function has redundant parameters and is related
to the classical Bessel function of the first kind \(\bessel_\alpha\) via
\begin{align}
    \nonumber
    \besselProj_n(\dims;x)
    \overset{\text{\citep[§22, Lem.~3]{mullerAnalysisSphericalSymmetries1998}}}&= \Gamma\Bigl(\frac{\dims}2\Bigr)\sum_{m=0}^\infty \frac{(-1)^m}{m!\Gamma(m+n+\frac{\dims}2)}\Bigl(\frac{x}{2}\Bigr)^{2m+n}
    \\
    \label{eq: standard bessel representation}
    \overset{\text{\citep[e.g.][9.1.10]{abramowitzHandbookMathematicalFunctions1964}}}&= 2^\lambda\Gamma\Bigl(\frac{\dims}2\Bigr) \frac{\bessel_{n+\lambda}(x)}{x^\lambda}.
\end{align}
That is the parameters \(n\) and \(\dims\) (or equivalently \(\lambda =\frac{\dims-2}{2}\)) may be merged
into the parameter \(\alpha\) of \(\bessel_\alpha\). 
Recall, by definition, Müller's Bessel function \(\besselProj\) is the projection of
\(e^{ir\cdot}\) onto the Gegenbauer polynomials (up to constants),
and \(e^{ir\cdot}\) is in \(L^2([-1,1], \weight)\). It therefore only requires
some care with the constants to verify the following representation of \(e^{it\cdot}\) \citep[see][§22,
Lem.~6]{mullerAnalysisSphericalSymmetries1998}:
\begin{equation}
    \label{eq: chebyshev expansion of fourier exponential}
    e^{ir t} = \sum_{n=0}^\infty i^n N_n^\dims \besselProj_n(\dims; r) \uPnorm_n(t).
\end{equation}
This identity is the key to the representation of the characteristic
function \[\Omega_\dims(t) \equalDef \E[e^{it\scp{\frac{v}{\norm{v}}}{U}}]\] 
of the uniform distribution
\(U\) on the sphere \(\sphere\). Recall that \(\Omega_\dims\) is used in the Schoenberg
characterization of \emph{stationary} isotropic kernels (Table~\ref{table:
characterizations}) and is therefore a natural starting point for the translation
into our characterization in Theorem~\ref{thm: stationary isotropic kernels}. In the proof of Theorem~\ref{thm: stationary isotropic kernels} we use the following lemma and the translation of Müller's Bessel function \(\besselProj_n\) into the classical Bessel function \(\bessel_\alpha\)
in \eqref{eq: standard bessel representation}.

\begin{lemma}[Characteristic function representation {\citep[§22, Lem.~2]{mullerAnalysisSphericalSymmetries1998}}]
    \label{lem: characteristic function of uniform distribution on sphere}
    For any \(x,y\in \real^\dims\) and \(s\in \real\) we have
    \begin{align}
        \Omega_\dims(s\norm{x-y})
        &= \besselProj_0(\dims; s\norm{x-y})
        \\
        &= \sum_{n=0}^\infty N_n^\dims \besselProj_n(\dims; s\norm{x})\besselProj_n(\dims; s\norm{y}) \uPnorm_n\bigl(\scp[\big]{\tfrac{x}{\norm{x}}}{\tfrac{y}{\norm{y}}}\bigr),
    \end{align}
    where the second equation is well defined for \(x=0\) or \(y=0\) as \(\besselProj_n(\dims; 0) = 0\)
    for \(n\ge 1\) and \(\besselProj_0(\dims; 0) = 1\) with \(\uPnorm_0\equiv 1\).
\end{lemma}
\begin{proof}
    The first identity follows from the Funk-Hecke formula (Lemma \ref{lem: Funk-Hecke}) and
    \[
    \begin{aligned}
        \Omega_\dims(s\norm{x-y})
        &= \E[e^{is\scp{x-y}{U}}]
        = \frac1{\vol{\sphere}}\int_{\sphere} e^{is\scp{u}{x-y}} d\sigma(u)
        \\
        &\overset{\text{Lem.~\ref{lem: Funk-Hecke}}} 
        = 
        \frac{\vol{\sphere[\dims-2]}}{\vol{\sphere}}\int_{-1}^1 e^{is\norm{x-y}t} (1-t^2)^{\frac{\dims-3}{2}} dt
        \\
        &= \besselProj_0(\dims; s\norm{x-y}),
    \end{aligned}
    \]
    where we used the definition of \(\besselProj_0\) and \(\uPnorm_0\equiv 1\)
    in the last step. The second identity is §22, Lemma 2 in
    \citep{mullerAnalysisSphericalSymmetries1998}. But since the proof is
    actually quite instructive, we provide it here. Observe that by the first
    line above we have
    \[
        \Omega_\dims(s\norm{x-y})
        = \frac1{\vol{\sphere}}\int_{\sphere} e^{is\norm{x}\scp{\frac{x}{\norm{x}}}{u}} \conj{e^{is\norm{y}\scp{\frac{y}{\norm{y}}}{u}}} \sigma(du).
    \]
    Consequently, using the Gegenbauer expansion of \(e^{ir\cdot}\) in \eqref{eq:
    chebyshev expansion of fourier exponential} with \(r= s\norm{x}\)
    and \(r= s\norm{y}\) respectively, we get
    \[\begin{aligned}
        &\Omega_\dims(s\norm{x-y})
        \\
        &= \!\!\sum_{n,m=0}^\infty\!\! i^{n-m}\besselProj_n(\dims; s\norm{x})\besselProj_m(\dims; s\norm{y})
        \scp[\Big]{
        \underbrace{N_n^\dims\uPnorm_n\bigl(\scp[\big]{\tfrac{x}{\norm{x}}}{\cdot}\bigr)}_{=\zonalK_n(\frac{x}{\norm{x}}, \cdot)}
        }{
        \underbrace{N_m^\dims \uPnorm_m\bigl(\scp[\big]{\tfrac{y}{\norm{y}}}{\cdot}\bigr)}_{=\zonalK_m(\frac{y}{\norm{y}}, \cdot) \quad \mathrlap{(\text{Thm.~\ref{thm: addition theorem}})}}
        }_{\sphere}
        \\
        &= \sum_{n=0}^\infty \besselProj_n(\dims; s\norm{x})\besselProj_n(\dims; s\norm{y}) \underbrace{\zonalK_n\bigl(\tfrac{x}{\norm{x}}, \tfrac{y}{\norm{y}}\bigr)}_{= N_n^\dims \uPnorm_n\mathrlap{\bigl(\scp[\big]{\tfrac{x}{\norm{x}}}{\tfrac{y}{\norm{y}}}\bigr)}}.
    \end{aligned}
    \]
    For the last equation we simply used the orthogonality of the spherical
    harmonics \(\sphHarmonic_n\) and \(\sphHarmonic_m\) for \(n\neq m\) and the
    fact that \(\zonalK_n\) is their reproducing kernel.
\end{proof}

\subsection{Spherical harmonics dimension}

For the discussion of the algorithmic complexity of Gaussian random function
simulation in Section~\ref{sec: GRF representation} we need the asymptotic
behavior of the dimension \(N_n^\dims\) of the space of spherical harmonics.
This is given in the following lemma.

\begin{lemma}[The dimension of spherical harmonics]
    \label{lem: spherical harmonics dimension}
    For \(\dims\ge 3\), the dimension of the space of spherical harmonics \(N_n^\dims=\dim \sphHarmonic_n\) is given by
\begin{equation}
    \label{eq: representations of spherical harmonics dimension}
    N_n^\dims
    = \frac{n+\lambda}{\lambda} \binom{n+2\lambda-1}{n}
    = \frac{n+\lambda}{\lambda}\uP_n(1),
    \qquad \lambda = \frac{\dims-2}{2},
\end{equation}
where \(\uP_n(1)\) is the value of the Gegenbauer polynomial at \(1\). And
we have \(N_0^2 =1\) and \(N_n^2 = 2\) for \(n\ge 1\).
For fixed dimension \(\dims \ge 2\) we have the asymptotic behavior
\[
    N_n^\dims \sim \tfrac{2}{(\dims-2)!}n^{\dims-2}
    \qquad \text{as } n\to \infty.
\]
\end{lemma}
\begin{proof}
    By Corollary 1.1.4 in \citep{daiApproximationTheoryHarmonic2013}
    \[
    N_n^\dims
    = \binom{n+\dims-1}{n} - \binom{n+\dims-3}{n-2}
    \]
    with the convention that \(\binom{n}{k} = 0\) for \(n < 0\). In particular, this
    implies \(N_0^2 = 1\) and \(N_n^2 = 2\) for \(n\ge 1\). For \(\dims \ge 3\) we have
    \[
        N_n^\dims
        = \frac{(2n+\dims-2)(n+\dims-3)!}{n!(\dims-2)!}
        = \frac{2n+\dims-2}{\dims-2}\binom{n+\dims-3}{n}.
    \]
The explicit representation follows from \(\lambda = \frac{\dims-2}{2}\) and the
formula for the Gegenbauer polynomials \citep[e.g.][§9.3]{mullerAnalysisSphericalSymmetries1998}.
For the asymptotic behavior we observe that the case \(\dims=2\) is trivial and for \(\dims \ge 3\) the equality
\[
    \frac{(n+\dims-3)!}{n!}
    = (n+\dims-3)\cdot \ldots \cdot (n+1)
    = n^{\dims-3} \prod_{k=1}^{\dims-3} (1+ \tfrac{k}{n})
\]
implies
\[
    \lim_{n\to \infty} \frac{\frac{(2n+\dims -2)(n+\dims - 3)!}{n!(\dims-2)!}}{\frac{2}{(\dims-2)!}n^{\dims-2}}
    = \lim_{n\to \infty} \frac{(2n+\dims-2)}{2n} \prod_{k=1}^{\dims-3} (1+ \tfrac{k}{n})
    =1
\]
and therefore
\[
    N_n^\dims \sim \tfrac{2}{(\dims-2)!} n^{\dims-2}
    \qquad \text{as } n\to \infty.
    \qedhere
\]
\end{proof}

\end{document}